\documentclass{amsart}
\usepackage[utf8]{inputenc}
\usepackage{enumerate}
\usepackage{amsmath,amsthm,amsfonts,amssymb,mathtools}
\usepackage{mathrsfs}
\usepackage{kotex}
\makeatletter
\@namedef{subjclassname@2020}{\textup{2020} Mathematics Subject Classification}
\makeatother

\usepackage{hyperref}
 \usepackage{enumitem}
\usepackage{capt-of}
\usepackage{relsize}

\newcommand{\Be}{\begin{equation}}
\newcommand{\Ee}{\end{equation}}
\newcommand{\Bea}{\begin{eqnarray}}
\newcommand{\Eea}{\end{eqnarray}}
\newcommand{\Bel}{\begin{align}}
\newcommand{\Eel}{\end{align}}
\newcommand{\Beas}{\begin{eqnarray*}}
\newcommand{\Eeas}{\end{eqnarray*}}
\newcommand{\Benu}{\begin{enumerate}}
\newcommand{\Eenu}{\end{enumerate}}
\newcommand{\Bi}{\begin{itemize}}
\newcommand{\Ei}{\end{itemize}}
\newcommand\supp{\operatorname{supp}}

\def\R{{\mathbb R}}
\def\Z{{\mathbb Z}}

\theoremstyle{plain}
\newtheorem{thm}{Theorem}[section]
\newtheorem{cor}[thm]{Corollary}
\newtheorem{lem}[thm]{Lemma}
\newtheorem{prop}[thm]{Proposition}

\theoremstyle{remark}
\newtheorem{rmk}[thm]{Remark}
\theoremstyle{definition}

\numberwithin{equation}{section}

\newcommand{\RNum}[1]{\uppercase\expandafter{\romannumeral #1\relax}}

\newcommand{\cA}{\mathcal A}
\newcommand{\B}{\mathcal B}

\usepackage{wrapfig}
\usepackage{tikz}

\newcommand{\cH}{\mathcal H}

\def\H{{\mathbf H}}
\newcommand{\re}{\operatorname{Re}} 
\newcommand{\im}{\operatorname{Im}} 

\newcommand{\by}{\mathbf y}
\newcommand{\rest}[1]{{\mathlarger |}_{#1}}

\newcommand{\epz}{\epsilon_\circ}
\newcommand{\sskip}[1]{\hspace{.#1pt}}

\newcommand{\epsc}{\epsilon_\circ}
\newcommand{\M}{\mathrm M}

\makeatletter
\@namedef{subjclassname@2020}{\textup{2020} Mathematics Subject Classification}
\makeatother

\title[Damping oscillatory integrals]{
Damping oscillatory Integrals of 
\\
Convex Analytic Functions}

\author[Sanghyuk Lee]{Sanghyuk Lee}
\author[Sewook Oh]{Sewook Oh}

\address{Department of Mathematical Sciences and RIM, Seoul National University, Seoul 08826, Republic of Korea}
\email{shklee@snu.ac.kr}

\address{Center for Mathematical Challenges, Korea Institute for Advanced Study, Seoul 02455, Republic of Korea}
\email{sewookoh@kias.re.kr}

\begin{document}

\begin{abstract}
    Let $\cH\subset \R^{d+1}$ be a compact, convex, analytic  hypersurface  of  finite type  with a smooth  measure $\sigma $ on $\cH$. Let $\kappa$ denote the Gaussian curvature on  $\cH$.    We consider the oscillatory integral  $(\kappa^{1/2}  \sigma)^\wedge$ with the damping factor  $\kappa^{1/2}$ and  prove the optimal decay estimate 
    \[ |(\kappa^{1/2}  \sigma )^\wedge(\xi)|\le C|\xi|^{-d/2}\] 
    for $d=2,3,$ and with an extra logarithmic factor  for $d=4$.  Our result provides an essentially complete answer, since  such decay estimates generally fail for $d \ge 5$, even for convex analytic hypersurfaces,  as shown by Cowling–Disney–Mauceri–Müller.
    Furthermore, we  prove  the same estimates for $(\kappa^{1/2+it} \sigma )^\wedge$ with $C$  growing polynomially in $|t|$.    As  consequences, we obtain the best possible estimates for the convolution, maximal, and adjoint restriction operators associated with $\cH$, incorporating the mitigating factors of  optimal orders. In particular,  for $d=2, 3$,  we prove the $L^2$--$L^{2(d+2)/(d+4)}$ restriction estimate with respect to the affine surface measure  $\kappa^{1/(d+2)} \sigma$.   This work was inspired by the  stationary set method due to Basu--Guo--Zhang--Zorin-Kranich.
\end{abstract}

\maketitle

\section{Introduction}

Let   $\cH \subset \mathbb{R}^{d+1} $  be  a smooth compact hypersurface and  $\sigma  $ be a smooth measure on $\cH$. 
The estimate for the Fourier transform  $\widehat{\sigma } $ is a fundamental subject in harmonic analysis.  The estimate  plays  a crucial role in various applications.   As widely understood among experts, it is essential for studying  problems associated with the measure $\sigma$ on $\cH$  such as  $L^2 $ restriction estimates,  $L^p $–$L^q $ convolution estimates, and  $L^p $ bounds on maximal averages. Additionally, the estimate  is important in the studies  of dispersive equations.

If  the Gaussian curvature is nonvanishing everywhere on $\cH$, the Fourier transform of ${\sigma } $ enjoys the best possible decay estimate: 
\[
|\widehat{\sigma }(\xi)|\le C |\xi|^{-d/2}
\]
for some constant $C$. 
This makes it possible to obtain the estimates on the optimal range for the convolution, maximal, and restriction operators defined by  $\sigma $ 
  (see, for example, \cite{Littman, Stein}). However, if the curvature vanishes at certain points, known as degeneracies, the decay estimate no longer holds and worsens. As a result, the best possible outcomes become not achievable. To address this issue, attempts have been made to recover the optimal decay by mitigating the problematic behavior near degeneracies through the  use of suitable damping factors (\cite{SS1, CDMM}), which  typically take the form of a power of the absolute value of the Gaussian curvature.

 To be specific, we assume that the surface  $\cH$ is given by the graph of a smooth function  $\phi:[-2,2]^d\to \mathbb R$, i.e., 
 \[ \cH=\{ (x, \phi(x)):  x\in [-2,2]^d\}. \]
Let $\mathrm{H} \phi(x)$ denote the Hessian matrix  $(\partial_{i}\partial_{j} \phi(x))_{1\le i, j\le d}$ of $\phi$. 
To simplify notation, we denote $\H\phi(x) = |\hspace{-1pt}\det \mathrm{H}\phi(x)|$. For $z\in \mathbb C$, we set 
\[
  \H^z \phi (x)= (\H\phi(x))^z.
\]
Instead of the Gaussian curvature,  we may use the determinant of the Hessian as a mitigating factor.   Indeed, 
 let $\kappa(x)$ denote the Gaussian curvature of $\cH$ at $(x, \phi(x))$. It is well known that  $|\kappa(x)|=b(x)\H\phi(x)$ for a nonzero smooth function $b$.   So, the damping factors  $|\kappa|^z$ and $ \H^z \phi$ play equivalent roles.

\subsection{Oscillatory integral with the damping factor}  Let  $\psi\in C_c^\infty(B_{1/16})$. Here $B_r\subset \mathbb R^d$ denotes the ball of radius $r$ centered at the origin.   We consider the measure  $\sigma^z  $ that is given by  
\[ (\sigma^z  , f)=\int  f(x,\phi(x))  \H^z \phi(x) \psi(x) dx.\]   
If the surface $\mathcal H$ has degeneracies,  the decay rate of the Fourier transform of $\sigma^z$  generally depends on the exponent $z$. 
There are various authors  who studied  
the estimate \Be 
\label{damped-opt}
|\widehat{\sigma^z  } (\xi)|\le C|\xi|^{-d/2}.
\Ee 
 
 In particular, the estimate \eqref{damped-opt} was first obtained by Sogge--Stein \cite{SS1} for  $ z = 2d  $, and later by Cowling--Mauceri \cite{CoM1} for  $ z = (d+3)/2  $ under the assumption that  $ \phi  $ is convex of finite (line)  type. A natural problem is to determine the \emph{optimal damping order}  $ \re z  $ for which the estimate \eqref{damped-opt} holds. 
From the perspective of the stationary phase estimate (e.g., \cite[Proposition 7.7.6]{H}), it appears reasonable to speculate that the mitigating factor  $ \mathbf {H}^{1/2} \phi(x)  $ is a natural candidate for the best possible decay order, namely  $ d/2  $. This corresponds to the estimate \eqref{damped-opt} when  $z = 1/2  $. On the other hand, the estimate \eqref{damped-opt} for $\re z=1/2$ with admissible growth in  $ |\im z|  $ has attracted particular interest. It is known that such an estimate  implies  the restriction, convolution, and maximal estimates, via interpolation along a suitable analytic family, for the operators defined by the surface  $ \mathcal{H}  $ with the damping factor  $ \mathbf{H}^\alpha \phi  $ of optimal orders  $ \alpha  $ (see Section~\ref{application} below).

In one dimension, the estimate  \eqref{damped-opt} for $\re z=1/2$  is less involved to prove  
when  $\phi$ satisfies a finite type condition (see \cite{KPV}). Uniform estimates over classes of phase functions also have been of interest. Notably, 
Oberlin \cite{Oberlin} obtained the estimate  with a uniform constant $C$ when the phase functions $\phi$ are polynomials of a fixed degree. 

 However, in higher dimensions, the problem becomes more intricate. As shown by Cowling--Disney--Mauceri--M\"uller  \cite{CDMM}, 
  the estimate \eqref{damped-opt} for $z=1/2$ is not generally valid. More precisely,  there are 
   finite type phase functions $\phi$ for which  the estimate  fails. It is not well understood yet how large the damping order should be to ensure that \eqref{damped-opt} holds.

Meanwhile, there also have been attempts to obtain the estimate  \eqref{damped-opt} for $z=1/2$  by restricting  the classes  of the phase $\phi$  (e.g., \cite{CoM, CDMM, CZ, CKZ, IM, BLN}).  In particular, Cowling et al. \cite{CDMM}  considered  the problem under the additional assumption  that $\phi$ is a smooth convex  function of finite type.  The sharp estimate  without the damping weight was earlier obtained  by Bruna--Nagel--Wainger  \cite{BNW}.  Even under the convexity assumption,  the finite type condition is necessary for the best possible decay. In fact,  \eqref{damped-opt} for $z=1/2$ does not hold for some flat convex hypersurfaces (e.g.,  see \cite{CZ}).   Cowling et al. \cite{CDMM}  proved that 
\[|\widehat{w \sigma }(\xi)|\lesssim |\xi|^{-d/2}\] 
if $0\le w\le \H^{1/2}\phi$ and $w$ is sufficiently smooth.    In view of the estimate \eqref{damped-opt}, the smoothness assumption on the damping factor $w$ appears  somewhat unnatural since  $\H^{1/2}\phi$ generally fails to be  smooth  on  the set $\{\H\phi=0\}$,  even if $\phi$ is convex and analytic. Consequently, this  makes the result in  \cite{CDMM} less practical in  applications. 
 
Nevertheless, it is still reasonable to attempt to obtain \eqref{damped-opt} for $z=1/2$ by further restricting  the class of phase functions, namely  to  convex analytic functions.  However,  when $d\ge 5$, the example provided in \cite[Example 1, p. 258]{CDMM}, which shows the decay order can be bigger than $2$,  rules out such a  possibility. When $2\le d\le 4$, the question of whether the estimate \eqref{damped-opt} for $z=1/2$ holds for convex analytic $\phi$ remains unsettled.  

The main purpose of this paper is to investigate these remaining cases $2\le d\le 4$.  For $d=2,3,$  we provide affirmative answers to the question. For $d=4$, we establish the estimate  with an additional logarithmic  factor  (see Theorem \ref{main00} below).

\subsection{Main result}   
The main difficulty in obtaining  $(1.1) $ for  $z = 1/2 $ stems from the lack of smoothness of  the damping factor  ${\H}^{1/2} \phi $,  which prevents us from applying integration by parts sufficiently many times near the singular set  $\{ \mathcal{H} \phi = 0 \} $. To overcome this, we adapt the stationary set method due to Basu–Guo–Zhang–Zorin-Kranich \cite{BGZZ}. 
Their key observation was that  if $S:\R\times\R^N\rightarrow\R$ is definable in an $o$-minimal expansion of $\R$,  then  \text{\it $S(\cdot, y)$  changes   monotonicity at most  $N$   times
with $N$} {\it being independent of $y$.}   This   was exploited to control the oscillatory integral by the measure of the stationary sets.

 As will be seen later in this paper, the stationary set method, which leverages monotonicity,  is less reliant on the smoothness of the amplitude function. 
 We utilize this advantage to obtain the desired estimate under the constraint of  limited smoothness by modifying the argument of Basu et al. On the other hand,  a convex function of finite type behaves as if it were a polynomial of mixed homogeneity (see, for example, \cite{Schulz}). After a suitable decomposition of the integral $\widehat{\sigma^z }$, this property  enables us to uniformly control  consequent  oscillatory integrals  via a normalization argument based on rescaling.  
 These two approaches constitute the key ingredients of our argument. By exploiting them, we obtain the following, which is the main result of this paper.

\begin{thm}\label{main00} 
Let $2\le d\le 4$, and $z=1/2+it$, $t\in \mathbb R$. 
  Suppose $\phi:[-2,2]^d\rightarrow\R$ is an analytic and convex function of finite type. Then,  for $|\xi|\ge 2$ we have 
\Be
\label{osc-est0}
|\widehat{\sigma^z  } (\xi)|\le C    (1+|t|)^{3}\begin{cases}    \  \     |\xi|^{-d/2},          \quad\quad  \  \  \quad &d=2,3, 
\\[4pt]
      (\log |\xi|) |\xi|^{-2},   \quad  \quad &d=4.     \end{cases}
     \Ee
\end{thm}

\subsection{Applications to the restriction, convolution, and  maximal estimates}
\label{application}
 Combined with Stein's interpolation along analytic families \cite[Theorem 4.1, Ch. V]{SW}), the inequality \eqref{osc-est0} has applications to the estimates for the (adjoint) restriction, convolution, and maximal operators.  When the surface  $\cH $ is degenerate, compared with the nondegenerate case,  the ranges of  $p $ and  $q $ for which   $L^p $– $L^q $ estimates hold shrink.  The problem of recovering the optimal ranges of  $p $ and  $q $, corresponding to the nondegenerate case,  using suitable damping factors has been studied by many authors (see \cite{SS1, KPV, Oberlin1, Oberlin2,   I1,  Shayya, CKZ,  Stovall,  IM, Li} and references therein).   In what follows, we address consequences of Theorem \ref{main00}  to those problems.

\subsubsection{$L^2$--restriction to the affine surface measure}
We first consider  the estimate 
\begin{equation}
\label{restriction0}
\| \widehat f\|_{L^2(\sigma^z  )} 
\le C \|f\|_p
\end{equation}
for all $f\in \mathcal S(\mathbb R^{d+1})$. 
If $\H\phi$ does not vanishes at any point in $\supp \psi$, the Stein--Tomas restriction theorem \cite{Stein, Tomas} tells that  the estimate \eqref{restriction0} for $z=0$ holds if and only if $p\in [1,p_\circ(d)]$, where 
\[ p_\circ(d):=\frac{2(d+2)}{d+4}.\]
 When $\phi$ is  a smooth convex function  of finite type,  the estimate  \eqref{restriction0}  for $z=0$ can be characterized by the decay of the Fourier transform of 
 $\sigma ^0$.  Iosevich  \cite{Iosevich} showed that  
  \eqref{restriction0} for $z=0$   holds for  $p\in [1, 2(r+1)/(r+2)]$ if and only if the estimate $|(\sigma ^0)^{\wedge}(\xi)|\le C |\xi|^{-r}$ holds for some $r>0$ (also, see \cite{BNW, greenleaf}).

When the Gaussian curvature of $\mathcal H$ vanishes, the estimate \eqref{restriction0} is no longer valid for all $p$ in $[1, p_\circ(d)]$. In such cases, the affine surface measure  has been considered to recover boundedness on the best possible range $[1, p_\circ(d)]$ (\cite{Oberlin2, CZ, Stovall}).  A related conjecture, known as the  \emph{$L^2$ affine restriction problem},  asserts  that \eqref{restriction0} holds for $z=1/(d+2)$ and $ p= p_\circ(d)$. The damping order  $1/(d+2)$  is the uniquely  optimal power.   

The estimate \eqref{damped-opt} with $\re z=1/2$ (and admissible constant $C$), 
combined with an $L^1$ estimate, implies  the estimate \eqref{restriction0} with $z=2/(p'd)$ for $1\le p\le p_\circ(d)$.\footnote{The order $2/(p'd)$  ensures affine invariance of the estimate \eqref{restriction0}.} This was observed by Kenig--Ponce--Vega \cite{KPV} in the case $d=1$, and their argument extends to higher dimensions $d\ge 2$ without modification;  see  \cite[Lemma 2.1]{CZ}. Combining Theorem \ref{main00} and this observation yields the following. 

\begin{cor}  Suppose  that $\phi$ and $\psi$ are given  as in Theorem \ref{main00}. If $2\le d\le 3$, the estimate 
 \eqref{restriction0} for $z=1/(d+2)$ holds for $p= p_\circ(d)$.  When $d=4$, \eqref{restriction0} with $z=2/(p'd)$ holds
  for $1\le p<p_\circ(d)$.
\end{cor} 

 When $d=2$,  there are numerous earlier results on the affine restriction problem.   Homogeneous polynomials $\phi$ were considered in \cite{CKZ}, while  the problem was studied for mixed homogeneous polynomials in \cite{Palle}.  
More recently, Li \cite{Li} studied  the problem for smooth $\phi$ in $\mathbb R^3$, and obtained \eqref{restriction0} in an almost sharp damping factor, i.e.,  for $z=1/4+\epsilon$ and $p= 4/3$. On the other hand, the problem with  hypersurfaces  with radial symmetry was studied  in \cite{Oberlin2, Stovall} (see, also, \cite{Shayya}).

\subsubsection{Convolution estimate}    In a similar vein,  one may also consider the  estimate    
\Be 
\label{conv}
  \|f\ast \sigma^z  \|_{q}\lesssim \|f\|_p. 
  \Ee
 When $\mathcal H$ has nonvanishing Gaussian curvature and the cutoff $\psi$ is  nontrivial and nonnegative, it is known (\cite{Littman})  that \eqref{conv} for $z=0$ holds if and only if 
 $(1/p, 1/q)$ is contained in the closed triangle $\mathcal T$ with vertices $(0,0)$, $(1,1)$, and   $((d+1)/(d+2), 1/(d+2))$ (see, for example,  \cite[Theorem 2]{Oberlin0}). 
 When the surface has degeneracy, it is easy to see that $\re z= 1/(d+2)$ is the optimal damping exponent for \eqref{conv} to hold.  
 A weak form of such a  estimate \eqref{conv} was studied in \cite{Oberlin1}.   
The following is a consequence of Theorem \ref{main00}. 

\begin{cor}\label{cor:conv} 
Let $d=2,3.$  Suppose that $\phi$ and $\psi$ are given as in Theorem \ref{main00}. Then, 
the estimate \eqref{conv}  with $ z=1/(d+2)$ holds if 
 $(1/p, 1/q)$ is in the closed triangle $\mathcal T$.  When $d=4$, \eqref{conv}  with $ z=1/(d+2)$  holds if $(1/p,1/q)$ is in the union of the interior of $\mathcal T$ and the line segment $\{(r,r)\in \R^2:0\le r\le1\}$. 
\end{cor}

\subsubsection{$L^p$ maximal bound}
We now consider  the maximal operator 
\[
\mathcal M^z f(x)=\sup_{t>0}|f\ast \sigma_{t}^z (x)|,
\]
where $\sigma_{t}^z$ denotes the $t$-dilation of $\sigma^z$, defined  by the relation $(\sigma_{t}^z, f)=(\sigma^z, f(t\,\cdot)).$  

The boundedness of the maximal operator  $M^0$, without the mitigating factor, has been studied  by many authors. 
For non-degenerate surfaces,  $M^0$ is bounded on $L^p$ if and only if $p>(d+1)/d$,  corresponding the celebrated results due to Stein \cite{St}  and Bourgain \cite{St}.   Also, there are various results for degenerate surfaces  (see, e.g., \cite{SS1, CoM, CoM1, IE, IE1}). In particular, for  the smooth convex surfaces, more satisfactory characterization of 
$L^p$ boundedness were obtained by Nagle--Seeger--Wainger \cite{NSW} and Iosevich--Sawyer--Seeger \cite{ISS}. 

 In degenerate cases, the admissible range of $p$ becomes narrower. As in the problems discussed above, it is natural to attempt to recover   maximal bounds  on the full range $p> (d+1)/d$ by introducing the mitigating factor  $\H^\rho \phi$ of a suitable order $\rho$.  When $d=1$,  Marletta \cite{Marletta} proved that $\mathcal M^{1/2}$ is bounded on $L^p$  for $p>2$, provided that the curve is degenerate but satisfies monotonicity and log-concavity. The damping order $1/2$ is optimal in the sense that the maximal operator $\mathcal M^\rho$ generally fails to be bounded on $L^p$  for all $p\in (2,\infty]$ if $\rho< 1/2$. That is to say,  when $\rho<1/2$,   there exists a curve of finite type  such that  $\mathcal M^\rho$ is not bounded on $L^p$ for some  $p\in (2,\infty]$.  

As an application of the estimate \eqref{osc-est0}, we obtain the following result, which provides the maximal bound with the damping factor $\H^\rho \phi$ of optimal order $\rho$ for $d=2,3,$ and $4$. 

\begin{cor}\label{cor:max} Let $2\le d\le 4$. Suppose that  $\phi$ and $\psi$ are as in Theorem \ref{main00}.
  Then, $\mathcal M^{1/(d+1)}$ is bounded on $L^p$ if and only if $p>(d+1)/d$.
\end{cor}

The damping order $1/(d+1)$ is also optimal. In other words, there exists a convex analytic  surface $\mathcal H$ for which $\mathcal M^{\rho}$ cannot  be bounded on $L^p$  for all $p>(d+1)/d$ if $\rho<1/(d+1)$ (see  Appendix \ref{appendix22}).

\subsubsection*{Organization of the paper} 
In Section~\ref{prelim}, we make several preliminary reductions for the proof of Theorem~\ref{main00}. In Section~\ref{sec:ssm},  we introduce a modification of the stationary set method to treat oscillatory integrals with weights, along with key notions such as o-minimal expansions and definable functions. In Section~\ref{sec:decomp}, using a decomposition adapted to convex functions of finite type, we reduce the proof of Theorem~\ref{main00} to establishing Proposition~\ref{lem:main}. Finally, Proposition~\ref{lem:main} is proved in Section~\ref{sec:main} using the modified stationary set method introduced earlier.

\subsection*{Notation}
For positive numbers $B, D$, we say $B\lesssim D$ if there exists a constant $C$ 
such that $B\le CD$. And  $B\sim D$ means that  $B\lesssim D$ and $D\lesssim B$.

\section{Preliminary reductions}\label{prelim}
In this section, to facilitate the proof of Theorem \ref{main00} we make some reductions. 
Writing $\xi=(\xi', \xi_{d+1})$,  to prove \eqref{osc-est0} we assume 
\[
|\xi_{d+1}|\ge c|\xi'|  
\] for some $c>0$.\footnote{For practical applications, this assumption is natural to impose, as rapid decay can otherwise be obtained through repeated integration by parts.}  The proof of 
\eqref{osc-est0} is much easier when $|\xi_{d+1}|< c|\xi'|$  (see Section \ref{main-proof}). 
  Instead of the Fourier transform of the measure $\sigma^z$, we consider
a modified  oscillatory integral 
\Be
\label{osi}
 \mathcal I^z(v,\lambda)= \int e^{i\lambda (\phi(x)+x\cdot v)}  \H^z \phi (x)\psi(x) dx.
  \Ee
The proof of Theorem \ref{main00} essentially reduces  to showing  the following.

 \begin{thm}\label{main} 
Let $2\le d\le 4$ and $\lambda\ge 2$.    Suppose that $\phi:[-2,2]^d\rightarrow\R$ is an analytic and convex function of finite type. Then, for a constant  $C$,  independent of $\lambda$ and $v$, we have 
\Be
\label{osc-est}
| \mathcal I^{1/2+ it}(v,\lambda)|  \le C     (1+|t|)^{3}\begin{cases}    \  \  \   |\lambda|^{-d/2},          \quad\quad  \  \  \quad &d=2,3, \\       {(\log|\lambda|)}|\lambda|^{-2},   \quad  \quad &d=4.     \end{cases}
     \Ee
  
\end{thm}

Before we proceed,  we show that the cut-off $\psi$ can be replaced by a semi-algebraic  function\footnote{A function is said to be semi-algebraic if its graph is given by a finite union of sets of the form $\{P=0\}\cap\bigcap_{j=1}^N \{Q_j>0\}$ for some polynomials $P,Q_1,\dots, Q_N$.} that is supported in a small neighborhood of a point $\omega\in B_{1/16}$.  This allows us to   apply the stationary set method, while  the oscillatory integral is properly localized near a point. 

\subsubsection*{Replacing  the cutoff  $\psi$ with a semi-algebraic function} It is not difficult to find a piecewise polynomial $C^2$ function $\psi_\circ$  such that $\psi_\circ\equiv 1$ on $B_{1}$ and $\supp \psi_\circ\subset B_{2}$.
Indeed, let  us define a  piecewise polynomial function $\eta_0$ on $\mathbb R$ by setting 
$\eta_0(r)=1$ for $|r|\le1$,  $\eta_0(r)=(2-|r|)^3(6|r|^2-9|r|+4)$ for $1<|r|<2$, and $\eta_0(r)=0$ for $|r|\ge 2$.  
Since $\eta_0(1)=1$, $\eta_0'(1)=\eta_0''(1)=0$, and $\eta_0(2)=\eta_0'(2)=\eta_0''(2)=0$, it is clear that $\eta_0\in C^2$. 
We now set 
\[\psi_\circ(x)=\eta_0(|x|), \] 
which is semi-algebraic,  and satisfies the desired properties.   

For a small $\epsilon_\circ>0$, let   $\{B(\omega_j,\epsc )\}_{j=1}^N$ be a collection of balls covering the support of $\psi$. Consider semi-algebraic functions 
\[\psi_j= \frac{\tilde \psi_j}{\sum_{\ell} \tilde\psi_\ell}, \quad 1\le j\le N,\]
forming a   partition of unity on  $\supp \psi$, where  $\tilde \psi_j=\psi_\circ((\cdot-\omega_j)/\epsc)$ and $\psi_j$ is defined to be zero whenever  $\sum_{\ell} \tilde\psi_\ell(x)=0$. 
Decomposing the integral $ \mathcal I^z$ with $\psi_j,$  
we have   $ \mathcal I^z= \sum_{j=1}^N \mathcal I^z_j$, where  
\[    \mathcal I^z_j(v,\lambda)= \sum_{j=1}^N \int e^{i\lambda (\phi(x)+x\cdot v)}  \H^z \phi(x) \psi_j\psi(x) dx. \] 
Since  $\psi\in C_c^\infty(B_{1/16})$, expanding $\psi$  in the Fourier series yields  
$
\psi(x)=\sum_{\mathbf k\in \Z^d}b_{\mathbf k}e^{ ix\cdot \mathbf k}
$
with $\sum_{\mathbf k\in \Z^d}|b_\mathbf k|\lesssim 1$. Thus, it follows 
\[
\mathcal I^z_j(v,\lambda)=\sum_{\mathbf k\in \Z^d}b_\mathbf k \int e^{i\lambda (\phi(x)+x\cdot (v+\lambda^{-1}\mathbf k))} \H^z \phi(x) \psi_j(x) dx.
\]
Obviously, we may also replace  $\phi$ by $\phi_\circ:=\phi\mathbf 1_{[-1,1]^d}$, which is a restricted analytic function, i.e., 
$\phi_\circ$ is  $\R_{an}$-definable. For the definitions of a restricted analytic function and an $\R_{an}$-definable function, we refer the reader forward to the beginning of Section \ref{sec:ssm}. Therefore,  Theorem \ref{main} follows once we prove the next theorem. 

\begin{thm}\label{main0}
   Let $2\le d\le 4$, and let  $\omega\in B_{1/16}$.  Suppose that $\phi$  is a restricted analytic and convex function of finite type and   $\psi \in \mathrm C^2_c(B(\omega, \epsilon_\circ))$ is a semi-algebraic function. 
   Then, \eqref{osc-est} holds
     with a constant $C$, independent of  $\lambda$, $\omega$, and $v$, if $\epsilon_\circ$ is small enough. 
\end{thm}

The above argument \emph{of replacing the cutoff function}  makes it possible  to {\it localize the problem} by decomposing the integral with definable cutoff functions. Since convexity is invariant under affine transformations, we may additionally  renormalize  the phase and amplitude functions so that they are well defined on a sufficiently large   ball centered at the origin (see {\it Remark} \ref{Note}).

\section{Stationary set method with weights}
\label{sec:ssm} 
In this section, we recall basic definitions,  such as  expansion of the real field, $o$-minimal structure, and their properties. 
Much of the material in this section overlaps with that in  \cite[Section 2]{BGZZ}, to which we refer the reader for details. We then prove a modified version of the stationary set method  to estimate oscillatory integrals with weights.

Consider the structure of the real ordered field $
\R=(\R,(0,1),(+,\cdot),\le)$, and its expansions:
\begin{align*}
&\R_{an}=\big(\R,(0,1),(+,\cdot\,, all \ restricted \ analytic \ functions),\le\big),\\
&\R_{an,exp}=\big (\R,(0,1),(+,\cdot\,, all \ restricted \ analytic \ functions,\exp),\le\big),
\end{align*}
 where $\exp:x\in \R\mapsto e^x\in \R$ is a function symbol. A function $f$ on  $\R^n$  is said to be a restricted analytic function if $f(x)=g(x)$ for $x\in [-1,1]^n$, where $g$ is an analytic function in a neighborhood of  $[-1,1]^n$,  and 
$f(x)=0$ otherwise.  A function $f:\R^{n}\rightarrow\R^m$ is definable in $\R_{an}$ ($\R_{an,exp}$, resp) if its graph is definable in  $\R_{an}$ ($\R_{an,exp}$, resp) (see, e.g.,  \cite[Section 2]{BGZZ} for the definition of  definability of a set in a structure).

\subsection{Definability  in $\R_{an}$}
We  now discuss some basic properties of definable functions in $\R_{an}$.  
Those will be useful when we try to verify definability of the phase and amplitude functions that arise later  through  decomposition and scaling.

\begin{prop}[{\cite[Remark 2]{Kurdyka}}]\label{prop:val} Let $f:\R^n\rightarrow\R^{n'}$ and $g:\R^{m}\rightarrow \R^{m'}$ be definable in $\R_{an}$. Then, 
    \begin{enumerate}[leftmargin=.8cm, labelsep=0.3 cm]
      \item[$(a)$]    $g\circ f$ is definable in $\R_{an}$ when $n'=m$. 
     \vspace{2pt}
        \item[$(b)$] $f\pm g$ and $fg$   are definable in $\R_{an}$ when $n=m$ and $n'=m'=1$.
    \end{enumerate}
\end{prop}

Consequently,  definability of functions is closed  under finite sums,  products, and composition. 
Note that semi-algebraic functions and restricted analytic functions are definable in $\R_{an}$. 
Combining these basic properties, we have the following lemma, which we use later. 

\begin{lem}\label{prop:definable}
Let  $f:\R^n\to\R$ be definable in $\R_{an}$. 
Then, the following functions are definable in $\R_{an}.$
\begin{enumerate}[leftmargin=.8cm, labelsep=0.3 cm]
\item[$(1)$] $\mathbf 1_{\{y \in \R^n:f(y)\in U\}}$ where $U$ is a finite union of intervals and points.
 \vspace{2pt}
\item[$(2)$] $F_{\le}(x,\beta):=\mathbf 1_{\{y\in \R^n: f(y)\le\beta\}}(x)$, \ $F_{<}(x,\beta):=\mathbf 1_{\{y\in \R^n: f(y)<\beta\}}(x)$.
 \vspace{2pt}
\item[$(3)$] $(x,T)\mapsto f(x,T):=f(Tx)$ where  $T$ is an affine transform from $\R^n$ to itself. 
 \vspace{2pt}
\item[$(4)$] $f^{-1} \mathbf 1_{\{f>0\}}$ and $f^{1/2}\mathbf 1_{\{f>0\}}$. Here,   by $g \mathbf 1_{E}$  we mean 
\Be\label{def}  g \mathbf 1_{E}(x)= 
\begin{cases} 
g(x),  &   x\in E,
\\[3pt]
 \  \  0,  &  x\notin E
\end{cases}
\Ee
 for $E\subset \mathbb R^d$ and $g:E\to \mathbb R$.
\end{enumerate}
\end{lem}

Throughout this paper, we frequently use the notation \eqref{def}. 
An affine map $T$ can be  identified with an element in 
$\R^{n^2}\times \R^n$. In fact, $Tx=Mx+v$ for some $(M,v)\in M_{n}(\R)\times \R^n$,  where $M_n(\R)$ is the collection of $n\times n$ matrices.  
For the rest of this paper, we regard an affine transformation $T$ as an element of  $M_{n}(\R)\times \R^n\simeq \R^{n^2+n}$, whenever we are concerned with definability.

\begin{proof}[Proof of Lemma \ref{prop:definable}]
 Note that $\mathbf 1_{\{x\in  \R^n:f(x)\in I\}}=\mathbf 1_I\circ f$.  Thanks to $(a)$ in Proposition ~\ref{prop:val}, it is enough  for  (1)  to show that $\mathbf 1_I:\R\rightarrow\R$ is definable in $\R_{an}$ if $I$ is a finite union of intervals and points. This is clear since the characteristic functions of a  finite union of intervals and points are semi-algebraic.
 
To show  (2), consider $\triangle_{\le}=\{(x_1,x_2)\in \R^2: x_1\le x_2\}$ and 
    $G(x,\beta)=(f(x),\beta)$. By $(a)$ in Proposition \ref{prop:val} we see that $F_{\le}=\mathbf 1_{\triangle_\le}\circ  G$  is definable in $\R_{an}$. 
   Similarly, one can show that $F_<$ is definable in $\R_{an}$.
   
For (3), we observe that the map $(x,M,v)\mapsto Mx+v$ is a polynomial map and hence definable in $\R_{an}$. Consequently,  by (a) in Proposition \ref{prop:val}  it follows  that $f(x,T)$ is definable in $\R_{an}$.
   
  To show  (4),  note that  the functions $g(r):= r^{-1} \chi_{\mathbb R_{>0}}(r) $ and $h(r):= r^{1/2}\chi_{\mathbb R_{>0}}(r) $ are semi-algebraic.  Thus, these functions are definable. 
  Since $g\circ f= f^{-1} \mathbf 1_{\{f>0\}}$ and $h\circ f=f^{1/2}\mathbf 1_{\{f>0\}}$, by $(a)$ in Proposition \ref{prop:val}  it follows that   $f^{-1} \mathbf 1_{\{f>0\}}$ and $f^{1/2}\mathbf 1_{\{f>0\}}$ are definable.    
\end{proof}

We conclude this subsection with a lemma to be used  in Section \ref{sec:main}, where we need to show definability
of derivatives of the phase and amplitude functions.  

\begin{lem}[\cite{Kurdyka}]\label{lem:kurdyka}
Let $f:\R^n\rightarrow\R$ be a definable function in $\R_{an}$, and let $U\subset \mathbb R^n$ be open and definable in $\R_{an}$. Suppose that $\partial_jf$ exists on $U$.  Then, $\partial_jf\mathbf 1_U$ is definable in $\R_{an}$.
\end{lem}

\subsection{$o$-minimal expansion} In this subsection we are  concerned with $o$-minimality, which played  a key role  in the work of Basu et al. \cite{BGZZ}.  Recall that 
an expansion of $\R$ is called  an $o$-minimal expansion of $\R$ if all definable subsets of  $\R$ of the structure are given by finite unions of points and intervals. From the perspective of oscillatory integral estimate,  the $o$-minimal property  is important in that it controls the number of monotonicity changes. More precisely,  if $f:\R\rightarrow\R$ is definable in an $o$-minimal expansion of $\R$, then there exist finitely many intervals $J_j$ such that $\R=\cup_j J_j$ and $f$ is either constant, strictly increasing, or strictly decreasing on each $J_j$(\cite[Theorem 1.2, Chapter 3]{Dries}). Furthermore, the following holds.
    
    \begin{prop} [{\cite[Proposition 2.8]{BGZZ}}]\label{monoton}
        Suppose that  $h:\R\times\R^n\rightarrow\R$ is definable in an $o$-minimal expansion of  $\R$.
       For $\by\in \mathbb R^n$, let $N(\by)$ denote the number of times  $h(\cdot,\by): \mathbb R\to \mathbb R$ changes monotonicity. Then $\sup_\by N(\by)<\infty$.
    \end{prop}
 
For the rest of the paper, we  consider only the aforementioned structures $\R_{an},$ $\R_{an,exp}$.
Our primary focus is on   $o$-minimality and  closedness  under integration in those structures. 
The following result is  due to van den Dries, Macintyre, and Marker \cite{DMM}.

\begin{thm} The structure 
$\R_{an,exp}$ is an $o$-minimal expansion of $\R$. 
\end{thm}

For our purpose, as to be seen later,  we need to handle 
  functions of the form
\[
I_f(y):=\int f(x,y)dx,
\]
while $f$ is definable in $\R_{an}$ and $f(\cdot, y)$ is  integrable in $x$. Unfortunately, definability of $f$  in $\R_{an}$  does not necessarily guarantee that of $I_f$. This was overcome in   \cite{BGZZ}  by considering 
the collection of constructible functions that was  introduced by 
Cluckers and Miller \cite{CM}. They showed that   the set   $ \{f: constructible\}$  is  closed  under  integration. 
Moreover, it is clear that 
    \[
    \{f:\text{definable in }\R_{an} \}\subset \{f:constructible\}\subset \{f:\text{definable in }\R_{an,exp}\}. 
    \]
Consequently, combining those facts, we obtain the following. 
  
    \begin{prop}\label{CMthm}
Let $f(x,y)$ be definable in $\R_{an}$ and $f(\cdot,y)$ be integrable for all $y$. Then
$I_f(y)$
is definable in $\R_{an,exp}$.
    \end{prop}
    
 Basu et al. \cite{BGZZ} focused on semi-algebraic functions since 
  they were concerned with uniform estimates for oscillatory integrals.    However,  thanks to Proposition \ref{CMthm},  we can apply Proposition \ref{monoton} for $h=I_f$ even when $f$ is merely  definable in $\R_{an}$. 
  
  The remaining of this section is devoted to proving a modification of the stationary set method \cite{BGZZ}. 
  For our purpose, we need to consider  the phase and amplitude functions definable in $\R_{an}$. Moreover, we allow the phase to include a logarithmic function related to the damping factor.
    
    \subsection{Stationary set method with weight} 
 In contrast to \cite{BGZZ}, where uniform estimates over semi-algebraic phases were obtained, we  consider a single  oscillatory integral given by  a specific analytic phase and do not pursue uniformity of the estimate over a class of phases.   
 
 To deal  with the damping factor with a complex exponent, we need to consider  the function $\mathbf 1_{\{\lambda\phi+t\log \H\phi\in[\beta,\beta+1]\}}$, which 
 is not generally constructible. This requires additional modification. In fact, we have the following.
 
\begin{thm}\label{complex SSM} Let $(x,\by)\in \mathbb R^n\times \mathbb R^m$. Let  $\Phi_1(x,\by),\Phi_2(x,\by)$, and $\mathbf a(x,\by)$ be  definable in $\R_{an}$.   Suppose that $\mathbf a(\cdot,\by)$ is integrable for all $\by$ and $\Phi_2>0$ whenever $\mathbf a\neq0$. For $L, \tau\in \mathbb R$, we set 
\[  I(\by, L,  \tau)=  \int e^{i(L \Phi_1(x,\by)+\tau\log \Phi_2(x,\by))} \mathbf a(x,\by)dx.\]
  Then, for $\tau\neq 0$, there is  a constant $C$, independent of  $\by,L,\tau$, such that
\Be
\label{main-est}
|I(\by, L,  \tau)|\le C(1+|\tau|) \int _0^\infty \sup_{\beta_1} S(\by, \beta_1,\beta_2,L, \tau') \frac{d\beta_2}{ \beta_2}
\Ee
where $ \tau'=e^{1/(\max(|\tau|,1))}$  and 
    \[
S(\by, \beta_1,\beta_2,L, s)=\int  (\mathbf 1_{\{(\bar x,\bar y):L \Phi_1(\bar x,\bar\by)\in [\beta_1,\beta_1+1]\}}\mathbf 1_{\{(\bar x,\bar y):\Phi_2(\bar x,\bar\by)\in [\beta_2, s\beta_2]\}} \mathbf a)(x,\by) dx. \]
\end{thm}

We remark that Theorem \ref{complex SSM} remains to be valid even if $\mathbf a$ is not supported on a compact set. However, we will not make use of this since our amplitude function  is assumed to have fixed compact support.

\begin{proof}
We handle  the cases $|\tau|\ge1$ and $0<|\tau|<1$, separately. First, we consider the case $|\tau|\ge1$. Note that
\begin{align*}
| I(\by,  L, \tau) |   &\sim \Big |\int_{-1}^0e^{i\beta_1}d\beta_1 \int_{-\tau/|\tau|}^0e^{i\beta_2}d\beta_2 \int e^{i(L\Phi_1+\tau\log \Phi_2)}\mathbf a\, dx\Big|.
\end{align*}
Here, we use the convention $\int_a^b=-\int_b^a$ and drop the variables $(x, \by)$ from $\Phi_1, \Phi_2$, and $\mathbf a$  for simplicity. Changing variables $(\beta_1,\beta_2)\mapsto (\beta_1-L\Phi_1,\tau\log(\beta_2/ \Phi_2))$, we see that
the right-hand side equals
\begin{align*}
   \Big|\int\int_{ e^{-1/|\tau|}\Phi_2}^{\Phi_2}\int_{L\Phi_1-1}^{L\Phi_1}\tau e^{i\beta_1}\beta_2^{-1+i\tau}\mathbf a\, d\beta_1d\beta_2dx\Big|.
\end{align*}
 Recalling the definition of $S(\by, \beta_1,\beta_2,L,\tau )$,  we have
\Be
\label{ineq:t>1}
  | I(\by, L,  \tau) | \lesssim |\tau|\int_{0}^\infty \Big  |\int e^{i\beta_1}S(\by, \beta_1,\beta_2,L, e^{1/|\tau|}) d\beta_1\Big| \frac{d\beta_2}{ \beta_2}.
\Ee

We next turn to the case  $0<|\tau|<1$. Observe that
\begin{align*}
   |I(\by, L, \tau)|
    &\sim \Big| \int_{-1}^0e^{i\beta_1}d\beta_1  \frac1{|\tau|}\int^{|\tau|}_0e^{i\beta_2}d\beta_2 \int e^{i(L\Phi_1+\tau\log \Phi_2)}\mathbf a\, dx\Big|. 
\end{align*}
Similarly, using the same change of variables as before, we have
\[
| I(\by, L,  \tau)|\lesssim \int_{0}^\infty \Big|\int e^{i\beta_1}S(\by, \beta_1,\beta_2,L,e) d\beta_1\Big|\frac{d\beta_2}{ \beta_2}.
\]
Therefore, by this and  \eqref{ineq:t>1},  the desired inequality \eqref{main-est} follows once we show  
\Be
\label{ss}
  \Big  |\int e^{i\beta_1}S(\by, \beta_1,\beta_2,L, \tilde \tau) d\beta_1\Big|\le C
    \sup_{\beta_1}|S(\by, \beta_1,\beta_2,L,\tilde \tau)|
\Ee 
with  $C$, independent of  $\by,\beta_2,L, \tilde \tau$. Clearly, we may assume $\tilde \tau>1$. 

To  this end, we first note that
\[
f(x, \by,\beta_1,\beta_2, L,  \tilde \tau)=(\mathbf 1_{\{ (\bar x,\bar\by):  L\Phi_1 (\bar x,\bar\by)\in  [\beta_1,\beta_1+1]\}} \mathbf 1_{\{(\bar x,\bar\by): \Phi_2(\bar x,\bar\by)\in  [\beta_2, \tilde \tau\beta_2]\}} \mathbf a)(x,\by)
\]
is definable in $\R_{an}$.  This can be shown using Proposition \ref{prop:definable} and Lemma \ref{prop:val}. 
Indeed, since  $\Phi_1$ and $\Phi_2$  are definable in $\R_{an}$,  by  $(2)$ (and $(3)$)  in Lemma \ref{prop:definable} it follows that  $h:=\mathbf 1_{\{(\bar x,\bar\by): L \Phi_1 (\bar x,\bar\by)\in [\beta_1,\beta_1+1]\}}$ and $g:=\mathbf 1_{\{ (\bar x,\bar\by): \Phi_2 (\bar x,\bar \by) \in [\beta_2, \tilde \tau\beta_2]\}}$ are  definable.
Thus,  $(b)$ in Proposition \ref{prop:val} shows that $f=h\hspace{.3pt}g\hspace{.5pt}\mathbf a $ is definable in $\R_{an}$. 

 Since $\mathbf a(\cdot, \by)$ is integrable,  so is $f(\cdot,\by,\beta_1,\beta_2,L, \tilde \tau)$ for all $\by, \beta_1,\beta_2,L, \tilde \tau$.  Thus, using  Proposition \ref{CMthm}  we see that 
 \[ S(\by, \beta_1,\beta_2,L,  \tilde \tau)=\int f(x, \by, \beta_1,\beta_2,L,  \tilde \tau) dx\] 
  is definable in $\R_{an,exp}$. 
 By $o$-minimality of $\R_{an,exp}$ and Proposition \ref{monoton}, there is an integer $N$, independent of $\by,\beta_2,L,\tilde \tau$,  such that $S(\by,\cdot, \beta_2,L,\tilde \tau)$ changes monotonicity at most $(N-1)$ times. So, for each $(\by,\beta_2,L,\tilde \tau)$,  there exists an integer $N'\le N$ and intervals $J_j=J_j(\by,\beta_2,L,\tilde \tau)$, $1\le j\le N'$,  such that 
 \[\R=\bigcup_{j=1}^{N'}J_j\] and $S(\by,\cdot,\beta_2,L,\tilde \tau)$ is monotonically increasing or decreasing on each $J_j$.  
Therefore, to show \eqref{ss} it suffices to prove that 
\Be
\label{ii}  
\mathcal I:=\Big|\int_{a}^b S(\beta)e^{i\beta}d\beta\Big|\le C\sup_{\beta\in (a,b)} |S(\beta)|\Ee
with $C$,  independent of $a,b,$ and $S$,  whenever  $S$ is either  increasing or decreasing on the interval $(a,b)$. We may assume that $S$ is monotonically increasing  since the other case can be handled by considering $-S$ instead of $S$. 

To show \eqref{ii},  write 
\[ [a,b]=\Big ( \bigcup_{m=0}^{M-1} [a+2\pi m, a+2\pi(m+1)]\Big) \cup [a+2\pi M, b] \] for some $M\ge 0$  such that the length  of the last interval $[a+2\pi M, b]$ is less than or equal to $2\pi$.  Here,  we regard the union $\bigcup_{m=0}^{M-1} $ as empty if $M=0$. 
Since $\int_s^{s+2\pi} e^{i\beta}d\beta=0$ for any $s\in \mathbb R$, we note that
\[ \mathcal J_m:=\int_{a+2\pi m}^{a+2\pi(m+1)}   S(\beta) e^{i\beta} d\beta = \int_{a+2\pi m}^{a+2\pi(m+1)}  \big( S(\beta) -S( a+2\pi m) \big) e^{i\beta}  d\beta. \]
  Consequently,  $|\mathcal J_m|\le 2\pi (  S(a+2\pi (m+1))-S(a+2\pi m) )$.  
Thefore, we have 
\[
\mathcal I\le  2\pi \sum_{m=0}^{M-1}  (  S(a+2\pi (m+1))-S(a+2\pi m) )  +    \int_{a+2\pi M}^b  |S(\beta)| d\beta.     
\]
Since $S$ is monotonically increasing, \eqref{ii} follows. 
\end{proof}

\section{Decomposition of the integral  $\mathcal I^{it +1/2}(v,\lambda)$}
\label{sec:decomp}

In this section,  we decompose the integral $\mathcal I^{it +1/2}(v,\lambda)$ in order to reduce the  proof of  the estimate \eqref{osc-est}    
to establishing  uniform estimates  for a family of oscillatory integrals whose phases and amplitudes belong to a certain class (see  Proposition \ref{prop:main} below). 
To this end, we renormalize the associated phase and amplitude functions via affine transformations provided  by John's lemma (see  Lemma \ref{john} below). This is where we heavily rely on  the \emph{finite type convexity} assumption. 
 Definability of the phase and amplitude functions along with their derivatives plays a crucial role in Section \ref{sec:main}, where we apply the stationary set method. 
 However, the discussion in this section relies solely on finite type convexity and does not make use of definability.
 
We begin by recalling some basic properties of the convex functions of finite type. Since $\phi=\phi_\circ\mathbf 1_{[-1,1]^d}$ where $\phi_\circ$ is a convex function of finite type on $[-2, 2]^d$, there exist an integer $k\ge 2$ and positive numbers $m$ and $M$ such that
\begin{align}
\label{phi-up} 
&|\partial^\alpha \phi(x)|\le M, 
\\
\label{phi-low}
\sum_{j=2}^k\frac{1}{j!}|&(\mathbf v\cdot\nabla)^j\phi(x)|\ge m
\end{align}
for all $x\in [-1, 1]^d$, $\mathbf v\in\{ v\in\R^d:|v|=1\}$, and $|\alpha| \le 2k+5$ (e.g.,  see \cite{BNW} and \cite{CDMM}).

The next lemma is a straightforward consequence of  \cite[Lemma 1.4]{CDMM} (also see Lemmas 3.1--3.3 in \cite{BNW}). 

\begin{lem}\label{lem-cdmm}
Let $f:[0,1]\rightarrow\R$ be a smooth convex function such that $f(0)=f'(0)=0$. Suppose that  
\Be
\label{finite1}
\sum_{j=2}^k \frac{1}{j!}|f^{(j)}(r)|\ge m, \qquad  \max_{0\le j\le 2k+5} |f^{(j)}(r)|\le M
\Ee
 for all $r\in [0,1]$ and some  $M, m>0$. 
Then,  the following hold with the implicit constants depending only on $k, m,$ and $M$$:$
\begin{align}
\nonumber 
f(r)&\sim \sum_{j=2}^k\frac{|f^{(j)}(0)|}{j!}r^j+ M r^{k+1}, 
\\
\label{CDMMlem2}
  f'(r)&\sim \sum_{j=1}^{k-1}\frac{|f^{(j+1)}(0)|}{j!}r^j+ M(k+1)r^{k}.
\end{align}
\end{lem}

Let $x, y\in B_{1/2}$, $x\neq y$, and $\mathbf v=(y-x)/|y-x|$.  Considering  $f(r)= \phi(r\mathbf v+x)-\phi(x)-\nabla\phi(x)\cdot \mathbf v r$, we apply  Lemma \ref{lem-cdmm}. Taking $r=|y-x|$,   by \eqref{CDMMlem2} we have 
\[
C|y-x|^{k-1} \le |\nabla\phi(y)-\nabla\phi(x)|   
\]
for some $C>0$ whenever  $x, y\in B_{1/2}$.  Consequently,  $x\to \nabla \phi(x)$ is injective and 
its inverse $(\nabla \phi)^{-1}: \nabla \phi(B_{1/2})\to B_{1/2}$ 
is H\"older continuous  of order $1/(k-1)$.

\begin{lem}\label{lemma:v}
   Let $-v\notin \nabla\phi(B_{1/4})$. Then, there is a positive constant $\delta$ 
   such that 
  \[  |\nabla\phi(x)+v|\ge \delta\]
   for  all $x\in B_{1/8}$.
\end{lem}

\begin{proof}[Proof of Lemma \ref{lemma:v}] 
Since $\nabla\phi(B_{1/4})$ is open and $\nabla\phi(\overline {B_{1/8}})$ is  a proper subset of $\nabla\phi(B_{1/4})$,  there clearly exists $\delta>0$ such that the $\delta$-neighborhood of $\nabla\phi(\overline {B_{1/8}})$  is contained in $\nabla\phi(B_{1/4})$. Consequently,  if $-v\notin \nabla\phi(B_{1/4})$, $|\nabla\phi(x)+v|\ge \delta$ for all $x\in B_{1/8}$. 
 \end{proof}

When $-v\notin \nabla\phi(B_{1/4})$,  Lemma \ref{lemma:v} implies that  the gradient of the phase $\phi(x)+x\cdot v$ of the integral $\mathcal I^{1/2+ it}(v,\lambda)$ is bounded away from zero. Thus,  the estimate  \eqref{osc-est}  becomes significantly easier to prove compared with the case when $-v\in \nabla\phi(B_{1/4})$. Indeed, 
  when $-v\notin \nabla\phi(B_{1/4})$,  the argument in Section \ref{sec:main} can be directly applied  without requiring sophisticated decomposition  (see {\it Remark} \ref{rmk:nonst} below).  Therefore, 
for the remainder of this paper, we assume 
\Be
\label{vcond}
-v\in \nabla\phi(B_{1/4}).
\Ee

\subsection{Normalization  of the convex function $\phi$}

Since $\phi$ is a convex function of finite type, 
as discussed above,  for $v$ satisfying \eqref{vcond} 
there exists a unique  $\omega=\omega_v\in B_{1/4}$ such that
\[
\nabla\phi(\omega)=-v.
\]
 
 We define 
\[  \Phi_v(x)=   \phi(x+\omega_v)+v\cdot x-\phi(\omega_v).\]
Note that $\Phi_v$  is also a convex analytic function of finite type on $B_{1/2}$.  Since $\Phi_{v}(0)=0$ and $\nabla\Phi_{v}(0)=0$,  there are constants $c_\phi$, $C_\phi$, depending only on $\phi$, such that 
\Be
\label{phiBphi}
 (c_\phi |x|)^{k} \le  \Phi_{v}(x) \le (C_\phi |x|)^2
\Ee
for $x\in B_{1/2}$.   Indeed, let $\theta\in \mathbb S^{d-1}$ and $f_\theta(2r)= \Phi_v(r\theta )$.  
Since $\Phi_v$ is also a convex function satisfying \eqref{phi-up} and \eqref{phi-low}, 
$f_\theta(r)$ satisfies \eqref{finite1} for some $m$ and $M$. By Lemma  \ref{lem-cdmm}, 
it follows that  $  c r^{k} f_\theta(1)  \le  f_\theta(r)\le C r^2   f_\theta(1) $ 
for some positive constants $c$ and $C$.  Since $f_\theta(1)\sim_{m, M} 1$ for all $\theta\in \mathbb S^{d-1}$, we  need only to take $r=2|x|$ and 
$\theta=x/|x|$ to get  \eqref{phiBphi}.

By \eqref{phiBphi}, we  have  $|x|\le c_\phi^{-1}h^{1/k}$ if  $ \Phi_{v}(x)\le h$. Let $0< h_\circ\le (4^{-1}c_\phi)^{k} $ be fixed. We take $ h_\circ$  to be small  enough (see {\it Remark} \ref{Note} below).  For $0< h<  h_\circ$, we set  
\[  \B_v^h= \{x\in  B_{1/2}:\Phi_{v}(x)< h/2 \}. 
\]
It should be noted that   $ \B_v^h$ is contained in  $ B_{1/4}$. It is known  \cite{BNW} that the collection of quasi-balls  $\{ \B_v^h\}$ possesses    
good geometric properties.  However, these will not to be employed in this paper.

Since $\Phi_{v}$ is convex,    the sub-level sets are convex sets with nonempty interior.  
Using the following lemma,  known as John's lemma  ({e.g., \cite[Lemma 4.1]{CDMM}}), 
 we can find an affine transformation that maps $\B_v^h$ to a set comparable to the unit ball $B_1$. 

\begin{lem}
\label{john}
Let $\Omega$ be a compact convex set with a nonempty interior in $\R^{d}$. Then there exists an invertible affine transformation $T$ such that
\[
B_1\subset T^{-1}\Omega\subset B_d.
\]
\end{lem}

Thus,  we have an invertible affine transformation  $T_v^h$ such that
\Be
\label{Tbv0}
B_1\subset  {\mathcal B}_{v,\circ}^h:=   (T_v^h)^{-1}  \B_v^h
\subset B_d.
\Ee
For our purposes, it is not necessary to know  the exact form of $T_v^h$.   However,  
since  $T_v^h(B_1)\subset \B_v^h \subset B_{c_\phi^{-1}h^{1/k}}$, it is clear that 
\Be
\label{Tb}
\|J(T_v^h)\|_{op}\le   C h^{1/k}, 
\Ee 
where $J(T)$ denotes  the Jacobian matrix of a map $T$.  We now  define an $h$ normalization of $\Phi_v$ by  
\Be
\label{phiv0}
\Phi_{v}^h(x)=h^{-1} \Phi_{v}( T_v^h x).
\Ee

\begin{rmk}\label{Note}  By taking $\epsilon_\circ$ in Theorem \ref{main0} to be sufficiently small,   $h_\circ$ can be made arbitrarily small. 
Furthermore,  making  $h_\circ$ to be small enough, we can ensure  that $\Phi_{v}^h$ is convex of finite type on a ball  $B_{R_\circ}$ of radius 
\[R_\circ\ge 100d\] as large as desired.  Thus, in what follows, we always assume that $\Phi_{v}^h$ is  convex of finite type on a  ball  centered at the origin with a sufficiently large radius $R_\circ$.
 \end{rmk}

\subsection{Properties of $\Phi_{v}^h$}   As observed in \cite{Schulz}, through an appropriate affine transformation, finite-type convex functions exhibit behavior similar to that of homogeneous polynomials with mixed homogeneity.  A key property of $\Phi_{v}^h$ is that these functions can be controlled in a uniform manner, independently of $v$ and $h$. For the remainder of this subsection, we establish several lemmas that enable us to obtain uniform control  over $\Phi_{v}^h$.

\subsubsection*{Upper bound on $\partial^\alpha_x \Phi_{v}^h$} 
The convex function $\Phi_{v}^h$ has uniformly bounded derivatives. This is basically a consequence of 
\eqref{phi-up}.

\begin{lem}\label{bounds}
Let $0<h<h_\circ$ and $|\alpha|\le 2k+5$. Then, 
for $x\in  B_{R_\circ}$, there exists  a constant $C>1$, independent of $h$ and $v$,  such that
\Be
\label{drv-bd}
|\partial_x^\alpha \Phi_{v}^h(x)|\le C.
\Ee
\end{lem}

\begin{proof} 
Let us write  $T_v^hy=\mathrm My+x_0$ where $\mathrm M=J(T_v^h)$ and $x_0\in \R^d$. 
Expanding  $\Phi_v(\cdot+x_0)$ into Taylor series up to degree $k$, we have 
\[
\Phi_{v}(y+x_0)=\sum_{|\alpha|\le k}a_\alpha y^\alpha
+R_{k+1}(y),
\]
where $R_{k+1}$ denotes the remainder. For simplicity, denote $P(y)=\sum_{|\alpha|\le k}a_\alpha y^\alpha$.
Let $y=\M x$ and define 
\[  P^h_v(x)=h^{-1}P(\M x), \quad R^h_v(x)=h^{-1}R_{k+1}(\M x). \] Then we have 
\[\Phi_{v}^h(x)= P^h_v(x)+  R^h_v(x).\]
From \eqref{phi-up} and \eqref{Tb}, it follows that $|\partial_x^\alpha R^h_v(x)|\lesssim h^{1/k}$   for $|x|\le R_\circ$.
So, there exists $h_\phi>0$, depending only on $\phi$, such that
\[
|\partial_x^\alpha R^h_v(x)|\le \frac 1{100d (k+1)!}
\] for all  $|x|\le {R_\circ}$ and 
$|\alpha|\le 2k+5$ provided that $0<h< h_\phi$. The bound \eqref{drv-bd}  is clear from \eqref{phi-up}  and \eqref{Tb} when   $h_\phi\le h<h_\circ$. Therefore, to show \eqref{drv-bd}, we may assume 
\[ 0<h<h_\phi.\]

Since the remainder $R^h_v$ already satisfies  the desired estimate, it is sufficient  to show that
\[ |\partial^\alpha P^h_v(x)|\lesssim 1 \] 
for $|x|\le {R_\circ}$.   Recall that $\Phi_v^h(x)\le 1/2$ for $x\in \mathcal B_{v,\circ}^h$ and $B_1\subset\mathcal B_{v,\circ}^h$. Thus, we have 
\Be
\label{upper}
 |P^h_v(x)|\le 1, \quad   x\in B_1. 
 \Ee
    We now  use an elementary fact about polynomials: If $p$ is a polynomial of degree at most $l$, i.e.,   $p(x)=\sum_{|\alpha|\le l} c_\alpha x^\alpha$, there is a constant $C_{d,l}$ such that
$\sum_{|\alpha|\le l} |c_\alpha| \le C_{d,l} \sup_{x\in B_1} |p(x)|.$   We write \[P^h_v(x)=\sum_{|\alpha|\le k} a_\alpha x^\alpha.\] 
Combining the above fact with \eqref{upper},   we have  $ \sum_{|\alpha|\le k} |a_\alpha|\lesssim_{d, k} 1.$  Therefore,  the desired inequality follows since $|x|\le R_\circ$. 
\end{proof}


\subsubsection*{Lower bound on  $|\nabla\Phi_{v}^h|$} 
For $0<h<h_\circ$, let us set 
\[
\mathcal D_{v,\circ}^h= \{x \in \mathbb R^d :1/2\le \Phi_{v}^h (x)\le 2\},
\]
which equals $ (T_v^h)^{-1}( \mathcal B_{v}^{4h}\setminus
 \mathcal B_{v}^{h})$. Using convexity and  \eqref{Tbv0}, one can show  that
\begin{equation}
\label{suppcond}
 \mathcal D_{v,\circ}^h\subset B_{9d}\setminus B_1.
\end{equation}
Indeed, note that $\mathcal D_{v,\circ}^h\cap B_1=\emptyset$ as is clear from the first inclusion in \eqref{Tbv0}.  Thus, \eqref{suppcond} follows if we show $\Phi_{v}^h (x)>2$ 
 whenever  $|x|\ge 9d$. 
 
 Let  $x_0$ denote the unique point such that $\Phi_{v}^h (x_0)=0$. By \eqref{Tbv0} it is clear that $x_0\in B_d$. 
Setting $\theta=(x-x_0)/|x-x_0|$, we consider the convex function 
\[  f(t)=\Phi_{v}^h(t\theta+ x_0). \] 
Let $t_1$ denote the unique positive number such that $t_1\theta+ x_0\in  S_d:=\{ x: |x|=d\}$. Thus, from  \eqref{Tbv0} \footnote{By \eqref{Tbv0} and \eqref{phiv0}, it follows that 
$ \Phi_{v}^h (x)\ge 1/2$   if $|x|\ge d.$} we have $f(t_1)\ge 1/2$, and $t_1< 2d$ since   $ x_0\in B_d$. Noting $f(0)=0$, we have  
\[   \frac{f(t_1)-f(0)} {t_1-0}=\frac{f(t_1)} {t_1}  >  \frac 1{4d}. \]
Denote by $t_2$ the number such that $t_2\theta+ x_0=x$. Since $|x-x_0|> 8d$, it is clear that $t_2> 8d$. On the other hand, thanks to convexity, it follows that  
\[ \frac{f(t_2)-f(0)} {t_2-0}\ge \frac{f(t_1)-f(0)} {t_1-0}.\]  Combining this and the above inequality yields 
$ f(t_2)>  t_2/{4d}.$
 Recalling $f(t_2)=\Phi_{v}^h(x)$ and $t_2> 8d$,   we conclude that  $\Phi_{v}^h (x)>2$  if $|x|\ge 9d$.

Using \eqref{suppcond}, we  obtain the following lower bound on $|\nabla\Phi_{v}^h|$, which plays an important role
in estimating the oscillatory integral. 

\begin{lem}\label{lem:lowbd}
Let $0<h<h_\circ$. Then, for all $x\in  \mathcal D_{v,\circ}^h$,  we have  
\[ |\nabla\Phi_{v}^h(x)|\ge 1/(20d).\]
\end{lem}

\begin{proof}   Let $x\in  \mathcal D_{v,\circ}^h$,  and let  $x_0$ denote, as  above, the point $x_0\in B_{d}$ such that $\Phi_v^h(x_0)=0$. 
Consider the convex function  
\[g(t)= \Phi_v^h(t(x-x_0)+x_0).\]
Since $g$ is convex,   $g(1)-g(0)\le g'(1)$. Thus, $(x-x_0)\cdot \nabla\Phi_v^h(x)\ge \Phi_v^h(x)$.
Consequently,  
\[
|\nabla\Phi_{v}^h(x)|\ge {|x-x_0|^{-1}} {\Phi_{v}^h(x)}.
\]
From \eqref{suppcond},  note that $|x-x_0|\le 10d$ if $x\in  \mathcal D_{v,\circ}^h$. Since 
 ${\Phi_{v}^h(x)}\ge 1/2$ for $x\in  \mathcal D_{v,\circ}^h$, the desired lower bound follows.    
\end{proof}

\subsubsection*{Derivatives of $\H  \Phi_v^h$}
By Lemma \ref{bounds},  the derivatives of $\H  \Phi_v^h$ are uniformly  bounded on the large ball $B_{R_\circ}$. Furthermore,  $\H  \Phi_v^h$ is nonnegative since $\Phi_v^h$ is convex. These facts can be utilized  to obtain an upper bound on the  gradient of  $\H^{1/2+ it}  \Phi_v^h$.  In particular,   it was shown by Glaeser \cite{Gl} that if a nonnegative function $f\in C^2$ on an open subset has its derivatives of order  $\leq 2 $ vanish at the zeros of  $ f  $, then  $ f^{1/2} \in C^1  $ (see also \cite{Dieudonne}).  Making use of the  arguments in \cite{Gl, Dieudonne}, we obtain the following.

\begin{lem}\label{root}     Let $0<h<h_\circ$ and $j\in\{1,\cdots, d\}$.   Then, we have
\Be
\label{root2}
    (\partial_{x_j} \H  \Phi_{v}^h(x))^2\le  C \H \Phi_{v}^h(x)
    \Ee
 for $x\in B_{R_\circ-1}$   with  a constant $C$, independent of $v,h$.   \end{lem}
  
This lemma shows that   the first-order derivatives of   $x\mapsto \H^{1/2+ it}  \Phi_{v}^h(x)$  are uniformly bounded outside the zero  set $Z_v^h=\{x\in B_{R_\circ}: \H  \Phi_v^h=0\}$.      Indeed, by Lemma \ref{root}  it follows that  there is a constant $C$, independent of $v,h$, such that 
    \[  | \nabla (\H^{1/2+ it}  \Phi_{v}^h)(x)|\le C(1+|t|) \] 
for $x\in B_{R_\circ-1}\setminus Z_v^h$.   

\begin{proof} 
    For simplicity, set $f(x)= \H\Phi_{v}^h(x)$. By Lemma \ref{bounds}, we have $|\partial_{x_j}f|,|\partial_{x_j}^2f|\le  M$ on $B_{R_\circ}$ for some $M>0$. Let $x\in B_{R_\circ-1}$ and  $r\in [-1,1]$. 
    Then,  by Taylor's theorem we get
    \[
    f(x_1,\cdots,x_j+r,\cdots,x_d)\le f(x)+\partial_{x_j}f(x)r+{M r^2}/{2}.
    \]
    Since $f$ is nonnegative,  taking $r=-\partial_{x_j}f(x)/M\in [-1,1]$ yields 
    the inequality $f(x)- (\partial_{x_j}f(x))^2/2M\ge 0.$
   Thus,   \eqref{root2} with $C=2M$ follows. 
\end{proof}

\subsection{Decomposition along dyadic  height $h$}  
Recall that $v$ satisfies \eqref{vcond}. By changing variables $x\to x+\omega_v$ in \eqref{osi},  we have 
\[ e^{-i\lambda \phi(\omega_v)}\mathcal I^{ 1/2+it}(v,\lambda)=   \int   e^{i\lambda \Phi_{v}(x)}  (\H^{1/2+it}\! \phi\hspace{1pt}\psi)(x+\omega_v) dx.\]

We now  choose a semi-algebraic cutoff function $\eta\in C^2(\R)$ such that both $\eta'$ and $\eta''$ are also semi-algebraic, $\supp\eta \in [-2,-1/2]\cup[1/2,2]$, and 
$
\sum_j\eta(2^{j} \cdot)\equiv 1. 
$
Indeed, it suffices to set  $\eta(r)=\eta_0(r)-\eta_0(2r)$, where $\eta_0$ is constructed in Section \ref{prelim}.

Let $\mathbb D=\{2^{j}: j\in \mathbb Z \}$. For  $h\in \mathbb D$, we set
  \[   \mathcal I^{h}_v(\lambda, t) =\int   e^{i\lambda \Phi_{v}(x)}  \eta(h^{-1}\Phi_{v}(x))  (\H^{1/2+it}\!\phi\hspace{1pt}\psi)(x+\omega_v) dx.\] 
 Since $\sum_{h\in \mathbb D} \eta(h^{-1}\Phi_{v}(x))=1$, 
we have 
\begin{align}
\label{mainobject}
| \mathcal I^{ 1/2+it}(v,\lambda)|\le \sum_{h\in \mathbb D} |  \mathcal I^{h}_v(\lambda, t)|.
    \end{align}
    
For $h\in \mathbb D\cap (0, h_\circ)$,  we define a phase function $\Psi_{\lambda,t}$ and an amplitude function  $\cA$, depending on $h,v$, 
by 
\begin{align}
\label{amp}
\Psi_{\lambda,t}(x)
&=\lambda h\Phi_{v}^h(x)+t\log \H \Phi_{v}^h(x),
\\
\label{psidefn}
\mathcal A(x)&=\eta(\Phi_{v}^h(x))\hspace{0.6pt} \H^{1/2}\hspace{-.4pt}\Phi_{v}^h(x)\hspace{0.6pt} \psi(T_v^hx+\omega_v).
\end{align}
Henceforth, we do not make  the dependence on $h,v$ explicit  for simpler notation.  However, it should be kept in mind that $\Psi_{\lambda,t}$ and $\cA$  always depend on $h,v$.

 Recall \eqref{phiv0} and  note  $
 \H \Phi_{v}^h(x)=
h^{-d} \det (J(T_v^h))^2\H \phi(T_v^h x+\omega_v)
$. Changing  variables $x\mapsto T_v^{h}x$ for the integral $ \mathcal I^{h}_v(\lambda, t) $, we have
\[  |\mathcal I^{h}_v(\lambda, t)| = h^{\frac d2} \big|\int e^{i  \Psi_{\lambda,t}(x)} \cA(x) dx\big|.\] 

The proof of the estimate \eqref{osc-est} is now reduced to showing the following. 

\begin{prop}\label{prop:main} 
Let $0< h <h_\circ$ and $1\le h |\lambda|$. If $h_\circ$ is small enough, then  there is a constant $C$, independent of $v,h$, $\lambda$,  and $t$, such that 
\begin{equation}
\label{eq:main}
   \Big |\int e^{i \Psi_{\lambda,t}(x)} \cA(x) dx\Big|
    \le C (1+|t|)^{3}(h\lambda)^{-2}.
\end{equation}
\end{prop}

Once we have Proposition \ref{prop:main}, the proof of \eqref{osc-est} is rather straightforward. 

\begin{proof}[Proof of Theorem \ref{main0}]
By \eqref{mainobject} and \eqref{eq:main} we have
\begin{align*}
 |\mathcal I^{ 1/2+it}(v,\lambda)| \lesssim \sum_{h\in \mathbb D,\, h< h_\circ} h^{d/2} \min \big\{1, (1+|t|)^{3}(h\lambda)^{-2}\big\}. 
\end{align*}
Consequently, dividing the sum into the cases $h\le |\lambda|^{-1}$ and $h>| \lambda|^{-1}$,  we have  the right-hand side  bounded by  a constant times 
\[    |\lambda|^{-d/2}+ (1+|t|)^{3}\int_{1/|\lambda|}^{2B_\phi} \frac{s^{d/2-1}}{ (\lambda s)^{2}}\, ds \]
for a positive constant $B_\phi$. 
Thus,  the estimate \eqref{osc-est} follows for $d=2,3$, and $4$.  
\end{proof}

To prove Proposition \ref{prop:main}, we make use of the lower bound on $\nabla\Phi_{v}^h$ in Lemma \ref{lem:lowbd}. To this end, we decompose the integral $ \mathcal I^{h}_v(\lambda, t)$ with suitable cutoff functions.  Let us set
\Be
\label{al-1} \cA_\ell (x)=   \cA (x) \frac{\mathcal(1-\eta_0) (20 d \sqrt d\, \partial_\ell \Phi_v^h(x))}{ \sum_{j=1}^{d} (1-\eta_0) (20 d \sqrt d\, \partial_j \Phi_v^h(x))} 
\Ee
for $\ell=1, \dots, d$. 
From this point forward, for any function expressed as a fraction, we assume the function equals zero whenever the denominator is zero.

Since  $1/2\le \Phi_{v}^h(x)\le2$ for $x\in\supp \cA$,  it is clear that  $\supp \cA \subset \mathcal D_{v,\circ}^h$.  Thus, 
by Lemma \ref{lem:lowbd} it follows  that   $(1-\eta_0) (20 d \sqrt d\, \partial_j \Phi_v^h(x))=1$ for some $j$ whenever $x\in \supp \cA$.  Consequently,  
\[
\textstyle\mathcal A (x)=\sum_{\ell=1}^{d}  \cA_\ell (x)  .
\]
By the support property of $\eta_0$, we have
\Be
\label{lbphase}
|\partial_{\ell} \Phi_{v}^h(x) |\ge 1/(20 d \sqrt d)
\Ee
for $x\in \supp \cA_\ell$. 
We further break $\mathcal A_l$ into two parts: one  away from,  and  the other near,  the set $Z_v^h$ where $ \H  \Phi_v^h$ vanishes.  For the former, the lower bound  \eqref{lbphase} will allow us to use integration by parts.

However, due to the term $t\log \H \Phi_{v}^h$ appearing in \eqref{amp}, the actual decomposition becomes slightly more involved. Note that 
\Be
\label{hoho}
\partial_\ell  \Psi_{\lambda,t}(x)=\lambda h\Big(\partial_\ell  \Phi_{v}^h(x)+\frac{t\partial_\ell  \H \Phi_{v}^h(x)}{ h\lambda \H\Phi_{v}^h(x)}\Big).
\Ee 
To utilize \eqref{lbphase}, we need an additional decomposition.  
Let us set 
\Be
\label{tilde-eta}
\tilde \eta (r)=\eta_0(r/(4C))-\eta_0( 80 d \sqrt d\, r),
\Ee
where $C$ is the constant appearing in Lemma \ref{bounds}.  
 We also set 
 \Be
\label{al}
    a_\ell (x)= \Big(1-\eta_0\big(h\lambda\H^{1/2} \Phi_{v}^h(x)\big)\Big) 
     \Big(1-\tilde \eta \Big(\frac{t\partial_{\ell} \H \Phi_{v}^h(x)}{h\lambda\H \Phi_{v}^h(x)}\Big)\Big).
\Ee

 We now break
\[
\cA_{\ell}=\cA_{\ell,0}+ \cA_{\ell,1}, 
\]
where 
\Be
\label{al-0}
\cA_{\ell,0}= \cA_{\ell}  (1-a_\ell),  \quad  \cA_{\ell,1}= \cA_{\ell} a_\ell.
\Ee
Note that $\tilde \eta (r)=1$ if 
 $|r| \in [1/(40 d \sqrt d),4C]$,  and  $ \tilde \eta(r) =0$ if $  |r| \not \in [1/(80 d \sqrt d),8C]$.   
 Thus,  
$|{t\partial_{\ell} \H \Phi_{v}^h}/{h\lambda\H \Phi_{v}^h|} \le  1/(40 d \sqrt d)$  or  $|{t\partial_{\ell} \H \Phi_{v}^h}/{h\lambda\H \Phi_{v}^h|} \ge 4C$ on  $\supp a_\ell$.  Using  \eqref{hoho} and \eqref{lbphase}, we have  
\Be
\label{low-psi}
|\partial_\ell  \Psi_{\lambda,t}|\gtrsim   \lambda h \min(1/(40 d \sqrt d), 3C )
\Ee 
on $\supp \cA_{\ell,1}$.  We now set 
\[ 
\mathcal I^{\lambda,t}_{\ell, \kappa}=\int e^{i  \Psi_{\lambda,t}(x)} \cA_{\ell, \kappa}(x)dx, \quad \kappa=0,1. 
\]
 
Now, Proposition \ref{prop:main} follows if we prove  the next proposition. 
\begin{prop}\label{lem:main} 
Let $0<h<h_\circ$, $\ell=1,\cdots,d$, and $\kappa=0,1$. If $h_\circ$  is small enough, 
then  there is a constant $C$, independent of $v,h$, $\lambda$,  and $t$, such that 
\begin{equation}
\label{stset}
    \big|  \mathcal I^{\lambda,t}_{\ell, \kappa}\big| \lesssim
   {(1+|t|)^{2+\kappa}}/{(h\lambda)^{2}}.
       \end{equation}
\end{prop}
In the next section, to prove  Proposition \ref{lem:main}, we use the fact that $\phi$ is a restricted analytic function and that $\psi$ is a $C^2$ semi-algebraic function.

\section{Proof of Proposition \ref{lem:main}}\label{sec:main}
 The oscillatory integrals  $\mathcal I^{\lambda,t}_{\ell, \kappa}$, $\kappa=0, 1$ in Proposition \ref{lem:main}  exhibit distinct differences. 
 On $\supp \cA_{\ell,0}$, $\H \Phi_v^h$ admits a good upper bound  (see \eqref{ineq:psi0} below)  but the amplitude function  $\cA_{\ell,0}$ has limited differentiability, so  integration by parts cannot be applied sufficiently many times.   On $\supp \cA_{\ell,1} $, by contrast,  $\H \Phi_v^h$ satisfies a lower bound (see \eqref{suppcavell1} below), which ensures 
 that  $\cA_{\ell,1}$ is smooth,  and  we also have the lower bound \eqref{low-psi}.  We handle both cases relying on the stationary set method. However, in the case $\kappa=1$, integration by parts  plays a significant role in obtaining the desired estimate. 
 
Since $t\mapsto  \mathcal I^{\lambda,t}_{\ell, \kappa}$ is continuous, we may assume $t\neq 0$.

\subsection{Proof of \eqref{stset} for $\kappa=0$}\label{sec:kappa0}  Recalling \eqref{al} and \eqref{al-0}, 
we note that 
\[
\H \Phi_{v}^h(x)\lesssim h^{-2}\lambda^{-2}, \  \text{ or } \   {\H \Phi_{v}^h(x)}\sim  |t|(h\lambda)^{-1}{|\partial_{\ell} \H \Phi_{v}^h(x)}| 
\]
for $x\in \supp\cA_{\ell,0} $. 
If the second holds, then applying \eqref{root2} with $j=\ell$ yields $\H \Phi_{v}^h\lesssim  |t|^2 /(h\lambda)^2$. Thus, $
\H \Phi_{v}^h\lesssim (1+|t|)^2 /(h\lambda)^2$  on $\supp\cA_{\ell,0} $.  As a result, we have  
\begin{equation}
\label{ineq:psi0}
|\cA_{\ell,0} (x)| \lesssim   \H^{1/2}   \Phi_v^h (x)  \lesssim    (1+|t|)(h\lambda)^{-1}
\end{equation}
for $x\in \supp\cA_{\ell,0} $. This inequality alone  is clearly  not sufficient to show the estimate \eqref{stset}.  To overcome  this, we make use of  Theorem \ref{complex SSM}.

For this purpose,  we need to verify $\R_{an}$-definability of  the phase and amplitude functions of the integral $\mathcal I^{\lambda, t}_{\ell, \kappa}$.  Rather than working directly with the original functions, we consider their extensions to a larger domain (i.e., with additional variables), of which definability is easier to show.
 We then apply Theorem \ref{complex SSM} to those extended functions.  After restricting them back to the relevant domain, we proceed to obtain the  desired estimate.

 To verify definability, we set $\mathbb I^d=[-1,1]^{d}$ and  consider 
 \[ \tilde\Phi (x, y, v,w):=  \begin{cases} 
\phi(x)+v\cdot y-\phi(\omega), &  (x, y, v, \omega) \in  Q
\\
\qquad \qquad  \   0,  &  (x,y, v, \omega) \not \in Q
\end{cases}, \] 
 where $Q=\mathbb I^d\times \mathbb I^d\times \mathbb I^d\times \mathbb I^d$.   Note that $\tilde \Phi$  is a restricted analytic function and its derivatives can also be defined as restricted analytic functions.
Putting 
\[  \by=(h, v, t,\lambda, \omega,T) \in \mathbb  \R \times \R^d\times \R\times \R \times  \R^d\times \R^{d^2+d}, \] 
where $T$ denotes an affine map identified as  an element in $\R^{d^2+d}$, we 
 define 
\[ 
\hat\Phi(x,\by):=  h^{-1}\mathrm 1_{(0,h_\circ)}(h)\sskip{7}\tilde\Phi (Tx+w, Tx, v,w).
\]
Consequently, for given $v$ and $h$, taking $\omega=\omega_v$ and $T=T_v^h$, we  have
\Be
\label{phiv}
\Phi_v^h (x)=\hat\Phi(x,\by)\rest{\omega=\omega_v, T=T_v^h}.
\Ee

Recalling \eqref{psidefn},  we also extend  $\cA$  to $\hat\cA$ by replacing  $\Phi_v^h(x)$ and $T_v^hx+\omega$ in \eqref{psidefn} by  $\hat\Phi(x,\by)$ and $Tx+\omega$, respectively.  Specifically, we set
\Be\label{aa2} \hat \cA(x,\by)= \eta(\hat\Phi (x,\by))\hspace{0.6pt}( \H^{1/2}\hspace{-.4pt}\hat\Phi\mathbf 1_{\{\H\hat\Phi>0\}})(x, \by)\hspace{0.6pt} \psi(Tx+\omega).\Ee 
The factor $\mathbf 1_{\{\H\hat\Phi>0\}}$ will not present any issues, as $\hat\Phi$ will be restricted to a set of convex functions.
Similarly to the above,  $\cA_\ell$ and $a_\ell$ are extended to $\hat\cA_\ell$ and $\hat a_\ell$, which are defined  by 
\begin{align}
\label{aaa2}
&\qquad  \hat\cA_\ell (x, \by)=  \hat \cA (x, \by) \frac{(1-\eta_0) (20 d \sqrt d\, \partial_\ell \hat \Phi (x,\by))}{ \sum_{j=1}^{d} (1-\eta_0) (20 d \sqrt d\, \partial_j  \hat \Phi(x, \by))},
\\[4pt]
\label{al2}
    \hat a_\ell &= \Big(1-\eta_0 \circ\big(h \lambda\, \H^{1/2} \hat \Phi \mathbf 1_{\{\H\hat\Phi>0\}}\big)\Big) 
     \Big(1-\tilde \eta \circ \big(\frac{t}{h\lambda} \frac{\partial_{\ell} \H  \hat \Phi }{\H \hat \Phi} \mathbf 1_{\{h\lambda\H\hat\Phi>0\}}\big)\Big)
\end{align}
({\it cf.}, \eqref{al}). 
  In the same manner, we extend 
 $\cA_{\ell, 0}$ to $\hat\cA_{\ell, 0}$ as follows: 
\Be
\label{al00}
\hat\cA_{\ell,0}(x,\by):= 
\mathbf 1_{\{\bar T:  \det J(\bar T) >0\}} (T) \hat \cA_{\ell}(x,\by)  (1-\hat a_\ell(x,\by)). 
\Ee
Consequently, we  also have 
\Be\label{al0} \cA_{\ell,0}(x)=\hat\cA_{\ell,0}(x,\by)\rest{\omega=\omega_v, T=T_v^h}.\Ee

Regarding definability of $\hat\Phi$ and $\hat\cA_{\ell,0}$, we have the following, which one can easily show using Proposition \ref{prop:val} and Lemma \ref{prop:definable}. 
 
\begin{lem}\label{definablity} The functions 
$ \hat\Phi$, $\H \hat\Phi= |\hspace{-1pt}\det \mathrm{H}_x \hat\Phi|$, and $\hat\cA_{\ell,0}$ are definable in $\R_{an}$. Here, 
$\mathrm{H}_x$ denotes the Hessian operator in $x$, i.e., $\mathrm{H}_x= (\partial_{x_j} \partial_{x_k})_{1\le j, k\le d}$.  
\end{lem}

\begin{proof}  
We first handle $\hat\Phi$.   By  $(1)$  in Lemma \ref{prop:definable} it is clear that $h^{-1}\mathrm 1_{(0,h_\circ)}(h)$ is definable. Similarly, 
$(x,\by) \mapsto 
\tilde\Phi (Tx+w, Tx, v,w)$, which is a composition of a restricted analytic function and a polynomial map, is definable. Consequently, by $(b)$ in Proposition \ref{prop:val}  it follows that $\hat\Phi$ is definable in $\R_{an}$. 

As for  $ \H \hat{\Phi}(x,\by) = |\hspace{-0.7pt}\det \mathrm{H}_x \hat{\Phi}(x, \by)|  $, since the map  $ s \mapsto |s|  $ is definable, (a) in Proposition \ref{prop:val} reduces the matter  to verifying that  $ (x, \by) \mapsto \det \mathrm{H}_x \hat{\Phi}(x, \by)  $ is definable. Since the determinant is a polynomial map, it suffices to show that  $ (x, \by) \mapsto \partial_x^\alpha \hat{\Phi}(x, \by)  $ is definable. Observe that
\[
\partial_x^\alpha \hat{\Phi}(x, \by) =  h^{-1}\mathrm 1_{(0,h_\circ)}(h)\sskip{7} \partial_x^\alpha ( \tilde \Phi(Tx + w, Tx, v, w)).
\]
Thus, we need only to show that $ (x,\by)\mapsto \partial_x^\alpha ( \tilde \Phi(Tx + w, Tx, v, w))  $  is definable. This is clear, since the map 
  is  a composition of  a restricted analytic function and a polynomial maps.

We now show definability of $\hat{\cA}_{\ell,0}$.  Recalling  \eqref{al00} and  \eqref{aaa2}, we need only to show that 
$\mathbf 1_{\{\bar T: \det J(\bar T)> 0\}}$, $\hat \cA$, $\hat \cA_\ell$, and $\hat a_\ell$ are definable in $\R_{an}$. First, it is clear that $\mathbf 1_{\{\bar T: \det J(\bar T)> 0\}}$ is definable since $\det: M_d(\R)\simeq \R^{d^2}\rightarrow \R$ is a polynomial map.  
To show definability of $\hat \cA$, recall \eqref{aa2}. Since $\H\hat \Phi$ is definable, by $(a)$ in Proposition  \ref{prop:val}  along with  (4) in Lemma \ref{prop:definable}   it follows that $\H^{1/2} \hat \Phi\mathbf 1_{\{\H\hat\Phi>0\}}$ is definable.  Since $\hat \Phi$, $\eta$, and $\psi$ are all definable,  making use of  $(a)$ in Proposition  \ref{prop:val}, we conclude that $\hat\cA$ is definable.

To show definability of $\hat \cA_\ell$, note that $\eta_0$ and $\partial_j \hat \Phi$ are definable. Hence,  so is $g_j:= (1-\eta_0) (20 d \sqrt d\, \partial_j  \hat \Phi)$.  
Setting $g=\sum_{j=1}^{d} g_j$, from \eqref{aaa2} we have  $\hat \cA_\ell=  \hat \cA g_\ell g^{-1} \mathbf 1_{g>0}$.   Consequently, by $(4)$  in Lemma \ref{prop:definable}, $\hat \cA_\ell$ is definable.  

Finally, we consider  $\hat a_\ell$. Recalling \eqref{al2}, we note that $\partial_{\ell} \H  \hat \Phi$ is definable, since  $\partial^\alpha_x\hat\Phi$ are definable.   
 As $\H\hat \Phi$ is definable,  by $(4)$  in Lemma \ref{prop:definable},    it follows  that  $(h\lambda{\H \hat \Phi})^{-1}{ \mathbf 1_{\{h\lambda\H\hat\Phi>0\}}}$ is definable. Thus, using $(a)$ in Proposition  \ref{prop:val}, we conclude that $\hat a_\ell$ is definable. 
\end{proof}

We first prove \eqref{stset} for $\kappa=0$ by applying  Theorem \ref{complex SSM} with the phase functions
\Be 
\label{phases}
  \Phi_1(x,\by)  =  \hat \Phi (x,\by),    \quad  \Phi_2(x,\by)  =  \H \hat \Phi (x,\by),
  \Ee
and the amplitude function  
\[ \mathbf a= \hat{\cA}_{\ell,0} (x,\by). \] 
Definability of these functions follows from Lemma \ref{definablity}. Moreover, integrability of  the function $\hat{\cA}_{\ell, 0}(\cdot,\by)$ is clear. Indeed, if $|J(T)|=0$, then $\hat{\cA}_{\ell, 0}\equiv 0.$  If $|J(T)|\neq0$,  then $\hat{\cA}_{\ell, 0}$ is bounded and supported in a  bounded set.  Therefore, Theorem \ref{complex SSM}  yields   the estimate \eqref{main-est} with  $\Phi_1$, $\Phi_2$, and $\mathbf a$ given as above. 
Recalling  \eqref{amp},  \eqref{phiv}, and \eqref{al0}, and taking $L=\lambda h$, $\tau=t$, $\omega=\omega_v$, and $T=T_v^h$ in the resulting estimate, we obtain
\[
|\mathcal I^{\lambda,t}_{\ell, 0}|=  \Big|\int e^{i  \Psi_{\lambda,t}(x)} \cA_{\ell, 0}(x)dx\Big|  \lesssim (1+|t|) \int _0^\infty  \mathfrak S_0(\beta_2)\frac{d\beta_2}{ \beta_2}, 
\]
 where 
\[  \mathfrak S_0(\beta_2)= \sup_{\beta}\Big|\int (\mathbf 1_{\{\bar x:\, \lambda h   \Phi_v^h (\bar x)  \in  [\beta, \beta+1]\}}\mathbf 1_{\{\bar x:\, \H  \Phi_v^h  (\bar x) \in [\beta_2, t'\beta_2]\}})(x) \cA_{\ell,0}(x) dx\Big|  \] 
and $t'=e^{1/\max(|t|,1)}$.

We note from \eqref{ineq:psi0} that  
\[ S(\beta_2):=\{x\in \supp \cA_{\ell, 0}:\, \H  \Phi_v^h (x) \in [\beta_2,t'\beta_2]\}=\emptyset \] if $\beta_2\ge C(1+|t|)^2/ (h\lambda)^{2}$ for a sufficiently  large constant $C$. Thus, 
\begin{align*}
\int_0^\infty  \mathfrak S_0(\beta_2)\, \frac{d\beta_2}{ \beta_2}
     \lesssim  \int_{0}^{C(1+|t|)^2/ (h\lambda)^{2}}  \widetilde{\mathfrak S}_0(\beta_2) \frac{d\beta_2}{\beta_2},
\end{align*}
where  
\[  \widetilde{\mathfrak S}_0(\beta_2)= \sup_{\beta}|\int (\mathbf 1_{\{ \bar x\in S(\beta_2):\, \lambda h   \Phi_v^h (\bar x)  \in  [\beta, \beta+1]\}}\cA_{\ell,0}(x) dx|.  \] 
By \eqref{ineq:psi0},  we have $|\cA_{\ell,0}(x) |\lesssim \beta_2^{1/2}$ for $x\in S(\beta_2)$.  Consequently, the desired estimate  \eqref{stset} for $\kappa=0$ follows via integration in $\beta_2$ 
provided that 
 \Be \label{a0} \sup_\beta |\{x\in \supp \cA_{\ell,0}: \lambda h\Phi_{v}^h(x)\in   [\beta, \beta+1]\}|\le C (h\lambda)^{-1} .  \Ee 
 Indeed, this estimate implies  $\widetilde{\mathfrak S}_0(\beta_2)\lesssim   \beta_2^{1/2} (h\lambda)^{-1} $. Hence, we obtain 
    \[  |\mathcal I^{\lambda,t}_{\ell, 0}|\lesssim  \frac{(1+|t|)}{\lambda h}\int_{0}^{C(1+|t|)^2/ (h\lambda)^{2}}  \frac{d\beta_2}{\beta_2^{1/2}}\lesssim  \frac{(1+|t|)^2}{(\lambda h)^{2}}.\]
 
In order to prove \eqref{a0}, we utilize  the following elementary lemma (see, e.g.,  \cite[Lemma 2]{CCW}). 

\begin{lem}\label{lem:set}  Let $I$ be a closed interval. Suppose $F:I\to \mathbb R$ satisfies $F'\ge 1$ on $I$. There is a constant $C$, independent of $I$ and $F$,   
such that  $|\{ s\in \mathrm I:  |F(s)|\le \rho  \}|\le C\rho$ for any $\rho>0$. 
\end{lem}

To show \eqref{a0}, we consider only the case $\ell=1$, as the other cases can be handled in the same manner.  Write $x=(s, x' )\in \mathbb R\times \mathbb R^{d-1}$.  Fixing $x'$, 
we set  
\[ S_{x'}= \big \{s\in \R:   \cA_{1,0}(s, x')\neq0  \text{ and } \lambda h\Phi_{v}^h(s,  x')\in   [\beta, \beta+1]\big \}. \]
For \eqref{a0}, it is sufficient to show that 
\[   |S_{x'}| \le C(\lambda h)^{-1}.\] 
Since both $ \hat \cA_{1,0}$  and $\hat \Phi$ are  definable,  Proposition \ref{monoton} implies that  $S_{x'}$ is a  union of at most $N$ points and 
intervals while $N$ independent of $x', v,$ $h$, and $\lambda$. 

We apply Lemma \ref{lem:set} to each of these intervals, taking  
\[F(s)= \Phi_{v}^h(s, x')- (h\lambda)^{-1}\beta.\] Recalling \eqref{lbphase},  note that  $|F'| \ge  1/(20 d \sqrt d)$ holds on all the intervals. Thus,  the desired inequality follows from Lemma \ref{lem:set}. \qed

\begin{rmk}  The argument above, in fact,  yields the estimate 
     \Be
     \label{ssss}
 \sup_{\beta}   \int_{\{\bar x\in \supp \cA_{\ell,0}:\, \lambda h  \Phi_v^h (\bar x)\in  [\beta, \beta+1], \, \H^{1/2} \Phi_v^h (\bar x)\le (h\lambda)^{-1}\}}\H^{1/2} \Phi_v^h (x)dx \lesssim (h\lambda)^{-2}
   .\Ee   
It is not difficult to see that  this estimate cannot  be improved  when  $\Phi_v^h$ is replaced  by  a general  convex function $\Phi$ of finite type.  
Indeed,  consider $ \Phi(x)=f(|x|)$ where 
\[ f(r)=\frac14 (r-1)^4 +r.\] 
Thanks to radial symmetry,  it is easy to see that the hessian matrix $\mathrm H \Phi (x)$ has eigenvalues $\lambda_1=3(|x|-1)^2$, $\lambda_2=\dots =\lambda_{d}=|x|^2-3|x|+3$. Thus, $\H \Phi (x)=3(|x|^2-3|x|+3)^{d-1}(|x|-1)^2$ and  $\Phi$ is convex 
on $B(e_1, \epz)$ for $\epz>0$.  However, 
\[  \int_{\{x\in B(e_1, \epz) :\lambda \Phi(x)\in  [\lambda, \lambda+1], \, \H^{1/2} \Phi (x)\le \lambda^{-1}\}}\H^{1/2} \Phi(x)dx \sim \lambda^{-2}\] 
if $\lambda>0$  is large enough. Thus,  it is not possible  to bound the left-hand side of \eqref{ssss}  above by $O(h\lambda)^{-d/2}$ when $d\ge 5$.  
\end{rmk}

\subsection{Proof of \eqref{stset} with $\kappa=1$}
\label{sec:nonpm}
Recalling \eqref{psidefn}, \eqref{al-1}, \eqref{al}, and \eqref{al-0}, we observe that $\H \Phi_{v}^h$ is bounded away from zero on $\supp \cA_{\ell,1} $. Thus,  $\cA_{\ell,1}$ is in $C^2$ 
since both  $\eta_0$ and $\tilde \eta$ are $C^2$ functions.   Integration by parts in $x_\ell$  gives 
\begin{align}
\label{eq:ibp}
    \int &e^{i \Psi_{\lambda,t}} \cA_{\ell,1} dx = -
    \int e^{i \Psi_{\lambda,t}} \mathcal D_\ell^2 ({\cA}_{\ell,1}) dx,
\end{align}
where    $ \partial_\ell=\partial_{x_\ell}$  and 
\[ \mathcal D_\ell  g =   \partial_\ell (({\partial_{\ell} \Psi_{\lambda,t} })^{-1} g) .\] 

To obtain \eqref{stset} with $\kappa=1$, we proceed as before, using  the stationary set method  for the right hand side of \eqref{eq:ibp} with the phases given by \eqref{phases}. However,  we must handle a different  amplitude function $\mathcal D_\ell^2 ({\cA}_{\ell,1})$.  Note that 
 \Be
 \label{deriv-a}
\mathcal D_\ell^2 ({\cA}_{\ell,1}) =\mathcal D_\ell^2 (1)\cA_{\ell,1}+ \mathcal D_\ell (1)\frac{3\partial_\ell\cA_{\ell,1}}{\partial_\ell\Psi_{\lambda,t}}  +\frac{\partial_\ell^2\cA_{\ell,1}}{(\partial_\ell\Psi_{\lambda,t})^2}.
\Ee
 
 By virtue of  support property of $\cA_{\ell,1}$,  $\mathcal D_\ell^2 ({\cA}_{\ell,1})$ is well defined without any ambiguity. However, additional care is required when we consider a definable extension as before.  Recalling \eqref{amp} and \eqref{psidefn}, we set 
\[ \hat \Psi (x,\by)=  \lambda h \hat\Phi(x,\by) + t  \log \H  \hat\Phi(x,\by), \]
and  \[ \hat\cA_{\ell,1}(x,\by)= 
\mathbf 1_{\{\bar T: \det J(\bar T) \neq0\}} (T) \hat \cA_{\ell}(x,\by) \hat a_\ell(x,\by), \]
where $\hat\cA_{\ell}$ and $\hat a_\ell$ are given by  \eqref{aaa2} and \eqref{al2}, respectively.  
Note that 
\Be \label{d-ell}
\begin{aligned}
\mathcal D_\ell (1)&=-\partial_\ell^2\Psi_{\lambda, t}/(\partial_\ell \Psi_{\lambda, t})^2, 
\\
\mathcal D_\ell^2 (1)&=-\partial_\ell^3\Psi_{\lambda, t}/(\partial_\ell \Psi_{\lambda, t})^3+ 3(\partial_\ell^2 \Psi_{\lambda, t})^2/(\partial_\ell \Psi_{\lambda, t})^4.
\end{aligned}
\Ee
This naturally leads us to consider \begin{align*}
\hat D_{\ell,1} (x,\by) &=-\frac{\partial_\ell^2\hat\Psi}{(\partial_\ell \hat\Psi)^2} (x,\by), 
\\[3pt]  \hat D_{\ell,2}(x,\by) &=-\frac{\partial_\ell^3\hat\Psi}{(\partial_\ell \hat\Psi)^3}(x,\by)+ 3\frac{(\partial_\ell^2 \hat\Psi)^2}{(\partial_\ell \hat\Psi)^4}(x,\by),
\end{align*}
provided that $\partial_\ell\hat\Psi(x,\by)$  exists and $\partial_\ell\hat\Psi(x,\by)\neq 0$. 
Clearly, $ \mathcal D_\ell (1)=\hat D_{\ell,1}$ and $ \mathcal D_\ell^2 (1)= \hat D_{\ell, 2}$ if we take  $\omega=\omega_v$ and $T=T_v^h$ for $\hat D_{\ell,1}$ and  $\hat D_{\ell, 2}$  whenever those are properly defined.   

Since $ \log \H  \hat\Phi$ is not defined on the set $\{\H  \hat\Phi\le 0\}$,  we introduce the characteristic function $\mathbf 1_{\{\H \hat \Phi>0\}}$.
To simplify notation, for a given function $G$, we denote 
\[ G_{\H}=G\mathbf 1_{\{\H \hat \Phi>0\}},\] 
where the right hand side is defined as in \eqref{def}.  We set 
\begin{align*}
(\hat D_{\ell,1})_{\H} &=-\frac{(\partial_\ell^2\hat\Psi)_{\H}}{(\partial_\ell \hat\Psi)^2_{\H}}, 
\\[3pt]
  (\hat D_{\ell,2})_{\H}&=-\frac{(\partial_\ell^3\hat\Psi)_{\H}}{(\partial_\ell \hat\Psi)^3_{\H}}+ 3\frac{(\partial_\ell^2 \hat\Psi)^2_{\H}}{(\partial_\ell \hat\Psi)^4_{\H}}. 
\end{align*}
 
 Recalling \eqref{deriv-a}, we now consider 
\Be
\label{am0}  \mathbf a=\Big( (\hat D_{\ell,2})_{\H} \hat \cA_{\ell,1} + (\hat D_{\ell,1})_{\H}\frac{3\partial_\ell\hat{\cA}_{\ell,1}}{(\partial_\ell\hat\Psi)_{\H}}+\frac{\partial_\ell^2\hat\cA_{\ell,1}}{(\partial_\ell\hat\Psi)^2_{\H}} \Big)\mathbf 1_{\{(\partial_\ell\hat\Psi)_\H\neq 0\}} .
\Ee
It is clear that  $\mathbf a=\mathcal D_\ell^2 ({\cA}_{\ell,1})  $ when we take  $\omega=\omega_v$ and $T=T_v^h$ for the right hand side, provided that $\partial_\ell\hat\Psi\neq 0$ and $\H\hat\Phi>0$ under the restriction.   

To apply Theorem \ref{complex SSM}, we now need to show that $\mathbf a$ is definable. One can show the definability of $\hat \cA_{\ell,1}$ in the same manner as that of $\hat \cA_{\ell,0}$, so we omit the details.   By Lemma \ref{lem:kurdyka}, it also follows that $\partial_\ell \hat \cA_{\ell,1}$ and $\partial_\ell^2 \hat \cA_{\ell,1}$ are definable.   
Since $\partial_\ell \H \hat\Phi$ and $\H \hat\Phi$ are definable, $(b)$ in Proposition  \ref{prop:val} and $(4)$  in Lemma \ref{prop:definable} show that 
\[ (\partial_\ell  \hat \Psi)_\H =  \lambda h \partial_\ell   \hat\Phi \mathbf 1_{\{\H \hat \Phi>0\}}  + t  \Big (   \frac{ \partial_\ell \H \hat\Phi}{ \H \hat\Phi}\Big)\mathbf 1_{\{\H \hat \Phi>0\}}\]
is definable.  Thus, $(\partial_\ell \hat\Psi)^m_{\H}$  is clearly definable for $m=2,3,4$.   Using  Lemma \ref{lem:kurdyka}, we also see that $(\partial_\ell^2 \hat\Psi)_{\H}$ and $(\partial_\ell^3\hat\Psi)_{\H}$ are definable.  Since $(\partial_\ell  \hat \Psi)_\H$ is definable,  combining these observations, we conclude that $\bf a$ is definable 
 thanks to the factor $\mathbf 1_{\{(\partial_\ell\hat\Psi)_\H\neq 0\}}$. In fact, one may use  $(4)$  in Lemma \ref{prop:definable} again.

We now apply Theorem \ref{complex SSM} with  $\Phi_1, \Phi_2$ given by \eqref{phases} and $\mathbf a$ given by \eqref{am0}. Consequently, by taking $L=\lambda h$, $\tau=t$, $\omega=\omega_v$, and $T=T_v^h$ in the resulting inequality as before, 
we obtain 
\[
|\mathcal I^{\lambda,t}_{\ell, 1}|=|\int e^{i \Psi_{\lambda,t}} \mathcal D_\ell^2 ({\cA}_{\ell,1}) dx|
\lesssim(1+|t|) \int _0^\infty \beta_2^{-1}\mathfrak S_1(\beta_2)d\beta_2,
\]
where 
\[
\mathfrak S_1(\beta_2):=\sup_{\beta_1}\int (\mathbf 1_{\{\bar x: \lambda h   \Phi_v^h (\bar x)  \in  [\beta, \beta+1]\}}\mathbf 1_{\{\bar x: \H  \Phi_v^h  (\bar x) \in [\beta_2, t'\beta_2]\}})(x)  |\mathcal D_\ell^2 ({\cA}_{\ell,1})(x)|   dx.
\]

To handle the weight part $|\mathcal D_\ell^2 ({\cA}_{\ell,1})| $, we use the next lemma. 

\begin{lem}\label{lem:psi1bound}
Let $0<h<h_\circ$. We have
    \Be
\label{psi1bound}
\begin{aligned}
| \mathcal D_\ell^2 ({\cA}_{\ell,1})|\lesssim (1+|t|)^2(\lambda h)^{-2}(\H \Phi_v^h)^{-1/2}.    
\end{aligned}
\Ee
\end{lem}

Assuming this for the moment, we proceed to prove \eqref{stset} for $\kappa=1$. From the definition of  $\mathfrak S_1(\beta_2)$, we may assume 
$\H  \Phi_v^h  \sim \beta_2$ since $\mathfrak S_1(\beta_2)=0$ otherwise. Thus, by  \eqref{psi1bound} and \eqref{a0},  it follows that 
\[
|\mathfrak S_1(\beta_2)|
\lesssim (1+|t|)^2(\lambda h)^{-3}\beta_2^{-1/2}.
\]
By \eqref{al},  we have $ \H  \Phi_v^h\ge ( \lambda h)^{-2}$ on the support of $ {\cA}_{\ell,1}$. Consequently, $ \beta_2\ge  ( \lambda h)^{-2}$ unless $\mathfrak S_1(\beta_2)=0$.  Therefore, we obtain 
\begin{align*}
   |\mathcal I^{\lambda,t}_{\ell, 1}|
    &\lesssim (1+|t|)^3\int_{(h\lambda)^{-2}}^{\infty}
   (\lambda h)^{-3}\beta_2^{-3/2}d\beta_2 \lesssim (1+|t|)^3(\lambda h)^{-2}.
\end{align*}

It remains to show  Lemma \ref{lem:psi1bound} to complete the proof of \eqref{stset} for $\kappa=1$.

\begin{proof}[Proof of Lemma \ref{lem:psi1bound}]
Clearly, it suffices to  show \eqref{psi1bound} on  $\supp \cA_{\ell,1}$.  Recall \eqref{deriv-a}. 
Thanks to \eqref{low-psi},  \eqref{psi1bound} follows if we show 
\begin{align}
\label{easy1}
| \mathcal D_\ell^m(1)| &\lesssim (1+|t|)^m(\lambda h)^{-m} (\H \Phi_v^h)^{- {m}/2},
\\[3pt]
\label{caderiv}
|\partial_{\ell}^m \cA_{\ell,1}|&\lesssim (1+|t|)^m(\H \Phi_v^h)^{-(m-1)/2} 
\end{align}
for $m=0, 1,2.$  Those estimates are trivial when $m=0$.  
 The first inequality \eqref{easy1} is easy to show.  Indeed, by \eqref{hoho} and a computation, we get
\begin{align*}
&\quad\partial_\ell^2\Psi_{\lambda,t}=\lambda h \partial_\ell^2 \Phi_v^h +t\frac{\partial_\ell^2\H \Phi_v^h}{\H \Phi_v^h}-t\frac{(\partial_\ell\H \Phi_v^h)^2}{(\H \Phi_v^h)^2},
\\[3pt]
\partial_\ell^3\Psi_{\lambda,t}&=\lambda h \partial_\ell^3 \Phi_v^h + t\frac{\partial_\ell^3\H \Phi_v^h}{\H \Phi_v^h}-3t \frac{\partial_\ell \H \Phi_v^h \partial_\ell^2 \H \Phi_v^h}{(\H \Phi_v^h)^2}+2t\frac{(\partial_\ell \H \Phi_v^h)^3}{(\H \Phi_v^h)^3}.
\end{align*}
Since $(1-\eta_0)(r)=0$ for $|r|\le 1$, we also note from  \eqref{al} that 
\Be
\label{suppcavell1}
(h\lambda)^{-1}\lesssim (\H \Phi_v^h)^{1/2}\lesssim 1
\Ee
on  $\supp \cA_{\ell,1}$.  
 By combining this,  \eqref{root2}, and  Lemma \ref{bounds},   we have 
$
|\partial_\ell^{m+1}\Psi_{\lambda,t}|\lesssim  (1+|t|)\lambda h (\H \Phi_v^h)^{-m/2}$ for $m= 1,2.$ 
Hence,  recalling \eqref{d-ell}, by  these inequalities and   \eqref{low-psi} we get  \eqref{easy1} for $m=1,2$. 

To show \eqref{caderiv}, we set 
\[  \mathcal A_\ast= 1-\eta_0\big(h\lambda (\H \Phi_v^h)^{1/2}\big) ,  \quad \mathcal A_\star= 1-\tilde \eta \big({t\partial_{\ell} \H \Phi_{v}^h}/({h\lambda\H \Phi_{v}^h})\big).\]
 Recalling \eqref{psidefn} and \eqref{al}, we  note from \eqref{Tb} and Lemma \ref{bounds}  that 
\[
\cA_{\ell,1}=\tilde \cA (\H \Phi_v^h)^{ 1/2}  
     \mathcal A_\ast  \mathcal A_\star
\]
for a smooth function  $\tilde \cA$ satisfying $|\partial^\alpha \tilde \cA|\lesssim 1$ for $|\alpha|\le 2$. Lemma \ref{bounds} and Lemma~\ref{root} show that  $|\partial_\ell^m (\H \Phi_v^h)^{ 1/2}  |\lesssim (\H \Phi_{v}^h)^{(1-m)/2}$ for $m=0,1,2$. Using this, we obtain  
\[
|\partial_\ell^m \mathcal A_\ast |\lesssim( \H \Phi_{v}^h(x))^{-m/2}, \quad m=0,1,2.
\]
Indeed, for $m=1,2,$ we use the fact that $h\lambda \sim(\H\Phi_{v}^h(x))^{-1/2}$ if  $\eta_0^{(k)}(h\lambda(\H \Phi_{v}^h(x))^{1/2})$ $\neq 0$, since $\eta^{(k)}$ is supported on $[-2,-1/2]\cup[1/2,2]$ for $k=1,2$.  Also, making use of  \eqref{suppcavell1}, one  can similarly  obtain 
\[    |\partial_\ell^m \cA_\star |\lesssim (1+|t|)^m ( \H \Phi_{v}^h(x))^{-m/2}, \quad m=0,1,2.\]
Combining all the estimates, we get \eqref{caderiv} for $m=1,2$. 
\end{proof}

Before closing this subsection,  we briefly explain how to handle the case when $-v\notin \nabla\phi(B_{1/4})$.

\begin{rmk}\label{rmk:nonst}  In this case, 
from Lemma \ref{lemma:v},  we have a good lower bound on the gradient of the phase, specifically, 
\[|\nabla\phi(x)+v|\ge c_\phi\]  
for  some constant  $c_\phi>0$ if $x\in B_{1/8}$. Consequently, the estimate \eqref{osc-est} becomes significantly easier to show.  
In fact, the following estimate holds:  
\Be
\label{strong}
 | \mathcal I^{1/2+ it}(v,\lambda)|  \le C     (1+|t|)^{3}\lambda^{-2}, \Ee
which is stronger than  the desired estimate  \eqref{osc-est}. 
Recalling \eqref{osi}, to show \eqref{strong} we follow the same lines of argument in Section \ref{sec:main} after replacing $ \Phi_v^h $, $\cA$ by \begin{equation*}
\begin{aligned}
 \tilde\Phi_v(x):=\phi(x)+v\cdot x,  \quad \tilde \cA(x)=   \H^{1/2} \tilde\Phi_v (x)\psi(x),
\end{aligned}
\end{equation*}
respectively.  It is clear that $ \H \tilde \Phi_v=\H\phi$ and 
$ |\nabla\tilde\Phi_v|\ge c_\phi$ on the support of  $\tilde \cA.$ Thus, the entire argument basically  reduces to the case $h\sim 1$ in Section \ref{sec:main}. 

Let $\tilde c=c_\phi^{-1}\sqrt d$.  To take advantage of  the lower bound  on $|\nabla  \tilde\Phi_v|$ and  apply Theorem \ref{complex SSM},  we proceed as before  by setting 
\begin{align}
\label{tilde-a}
\tilde\cA_{\ell}&=    \tilde\cA\, \frac{(1-\eta_0) (\tilde c\, \partial_\ell  \tilde \Phi_v)}{ \sum_{j=1}^{d} (1-\eta_0) (\tilde c\, \partial_j  \tilde \Phi_v)},  
\\
\label{tilde-aa}
   \tilde a_\ell &= \Big(1-\eta_0\big(\lambda \H^{1/2}\phi \big)\Big) 
     \Big(1-\tilde \eta \Big(\frac{t\partial_{\ell} \H \phi}{\lambda\H \phi}\Big)\Big)
     \end{align}
     for $\ell=1, \dots, d$ ({\it cf.}, \eqref{al-1} and \eqref{al}). Here, we choose $\tilde \eta$ by adjusting the associated constants in \eqref{tilde-eta} so that 
$|{t\partial_{\ell} \H \phi}/{\lambda\H \phi|} \le c_\phi/2\sqrt d$  or  $ |{t\partial_{\ell} \H \phi}/{\lambda\H \phi|} \ge 2 \sup\{x\in\supp \psi:   |\nabla\tilde\Phi_v(x)|\}  $ 
on $\supp  \tilde a_\ell$.   We also set 
\[
\tilde\cA_{\ell,0}  := \tilde \cA_{\ell}  (1-\tilde a_\ell), \quad \tilde\cA_{\ell,1}:= \tilde \cA_{\ell} \tilde a_\ell.
\]
Thus, the matter reduces to proving  the estimates
\[
|\int e^{i(\lambda\tilde\Phi_v+t\log \H\tilde\Phi_v)}\tilde\cA_{\ell,\kappa}dx|\lesssim (1+|t|)^{2+\kappa}\lambda^{-2}, \quad  \kappa=0,1.
\]

One can establish these estimates by routinely following the lines of argument in Section \ref{sec:kappa0} and \ref{sec:nonpm},  in particular  by applying Theorem \ref{complex SSM} with 
\[ \Phi_1= \tilde \Phi_v, \  \Phi_2=\H\phi, \text{  and } \mathbf a=\tilde\cA_{\ell, \kappa}.\]  Since $|\partial_\ell \tilde \Phi_v (x)|\ge c_\phi$ on $\supp \tilde\cA_{\ell, 0}$,  the estimate for $\kappa =0$  readily follows by repeating the argument in Section \ref{sec:kappa0}. 
On the other hand, note that the phase function $  \tilde \Psi_{\lambda,t}:=\lambda \tilde \Phi_v+ t\log \H\phi$ satisfies 
\[ |\partial_\ell  \tilde \Psi_{\lambda,t}(x)|\gtrsim  \lambda
\]
 on $\supp \tilde \cA_{\ell, 1}$ ({\it cf.}, \eqref{low-psi}). Thus, by applying integration by parts and invoking the argument from Section \ref{sec:nonpm}, one obtains the desired estimate for $\kappa = 1$.
\end{rmk}

\subsection{Proof of Theorem \ref{main00} when $|\xi_{d+1}|< c|\xi'|$} \label{main-proof}
Finally, to complete  the proof  of Theorem \ref{main00},  we consider  the case 
\[|\xi_{d+1}|< c|\xi'|,\]
which is not covered by Theorem \ref{main}.  This can be easily handled in the same manner as in {\it Remark} \ref{rmk:nonst}. 
However, for completeness, we include some details.

By the argument in Section \ref{prelim},   we may assume that the cutoff $\psi$ is semi-algebraic. 
Setting $\lambda=-|\xi'|$,  $u=-\xi'/|\xi'|$, and $c_0=  -\xi_{d+1} /|\xi'|$ gives  
\[  \widehat{\,\,\sigma ^{1/2+it}} (\xi)=   \int e^{i(\lambda \bar \phi(x)+t\log \H\phi(x))} \bar \cA(x) dx, \]
where  
\[ \bar \phi(x) =  u \cdot  x +  c_0 \phi(x), \quad \bar\cA= (\H^{1/2}\phi)\psi.\] 
Thus, taking $c$ sufficiently small,  we have  
\[  |\nabla \bar \phi(x)|\ge 1/2,   \quad  x\in B_{1/16}.  \]

By virtue of this lower bound, the remainder of the argument  is almost identical to that in {\it Remark} \ref{rmk:nonst}.  Indeed, there is an $\ell$ such that  $|\partial_\ell \bar \phi(x)|\ge (2\sqrt d)^{-1}$ for $x\in B_{1/16}$.  Let $\bar c=2\sqrt d$. 
Define $\bar \cA_{\ell}$ as in \eqref{tilde-a}, by replacing $\tilde c$ and $\tilde \Phi_v$ with $\bar c$ and $\bar \phi$ respectively.  
We also define 
 $   \bar a_\ell =\tilde a_\ell $ (given by \eqref{tilde-aa}) analogously. Here,  as before, $\tilde \eta$ is chosen  by adjusting constants in \eqref{tilde-eta} so  that 
$|{t\partial_{\ell} \H \phi}/{\lambda\H \phi|} \le \bar c/2$   or  $ |{t\partial_{\ell} \H \phi}/{\lambda\H \phi|} \ge 2 \sup\{x\in\supp \psi:   |\nabla\bar \phi(x)|\}  $ 
on $\supp  \bar a_\ell$. We also set 
\[ 
\bar \cA_{\ell, 0}= \bar \cA (1-\bar a_\ell),   \quad \bar \cA_{\ell, 1}= \bar \cA \bar a_\ell.
\]
Therefore, it is sufficient show that
\[
|\int e^{i(\lambda \bar \phi(x)+t\log \H\phi(x))}\bar\cA_{\ell, \kappa} dx|\lesssim (1+|t|)^{2+\kappa}\lambda^{-2}, \quad  \kappa=0,1. 
\]

 Now, as in  {\it Remark} \ref{rmk:nonst},   a routine adaptation of  the arguments in Section \ref{sec:kappa0} and \ref{sec:nonpm} yields those estimates  by applying Theorem \ref{complex SSM} with $\Phi_1=\bar \phi$, $\Phi_2=\H\phi$, and $\mathbf a=\bar\cA_{\ell, \kappa}$.  We omit the details.  \qed

\appendix

\section{Maximal estimates}

For the reader's convenience, we provide the proofs of Corollary \ref{cor:conv} and \ref{cor:max}. 
 
\subsection{Proof of Corollary \ref{cor:conv}}

When the cutoff $\psi$ is  nontrivial and nonnegative, the failure of the $L^p$--$L^q$ bound on the convolution operator 
\[ T: f\mapsto f\ast\sigma^{1/(d+2)}\] 
 for $(1/p, 1/q)\not\in \mathcal T$ follows from that of  the convolution operator 
defined by a surface with nonvanishing curvature (e.g., \cite[Theorem 2]{Oberlin0}). Indeed, since the cutoff $\psi$ is  nontrivial and nonnegative, there is a neighborhood $V\subset \supp \psi$ such that $\psi>0$ on $V$.    By finite type convexity of  $\phi$,  there is a point $y_0\in V$ such that $\H\phi\ge c$ on a small neighborhood $V_0$ of $y_0$ for some positive $c$. This can easily be seen using the asymptotic expansion due to Schulz \cite{Schulz}. Let $\psi_0$ be a nonnegative function such that $0\le  \psi_0\le1$,  $\psi(y_0)=1$, and $\supp \psi_0\subset  V_0$. Set $\sigma_0 = \psi_0 \sigma^{0}$, which is a measure on a surface with nonvanishing curvature.    Since $\sigma_0\le   c^{-1/(d+2)} \sigma^{1/(d+2)},$ the  failure of the $L^p$--$L^q$ on  $f\mapsto f\ast\sigma_0$ implies that of $T.$

\subsubsection*{Proof  Corollary \ref{cor:conv}} 
We first consider the cases $d=2,3$.
Let $\beta\in C_c^\infty((1/2,2))$ such that $\sum_{j=-\infty}^\infty \beta(2^{-j}\cdot) =1$ on $\mathbb R_{>0}$. Let  $P_j$ denote the operator 
given by 
\[\widehat {P_j f}=  \beta(2^{-j}|\cdot|)\widehat f.\] We also denote by $K_j$ the kernel of $P_j$, i.e., $K_j=\beta(2^{-j}|\cdot|)^\vee$.

Since $\sigma^{{1}/(d+2)}$ is a finite measure, $\|f\ast \sigma^{{1}/(d+2)} \|_{p} \le C   \|f\|_{p}$ for $1\le p\le \infty$. By interpolation  
it suffices to show  that $T$ is bounded from $L^{d+2}$ to $L^{(d+2)/(d+1)}$. To this end, we only need to show 
\Be \label{lp} \| P_jf\ast \sigma^{{1}/(d+2)} \|_{d+2} \le C   \|f\|_{\frac{d+2}{d+1}}\Ee
for $j\ge 0$. 
Since $\frac{d+2}{d+1}< 2< d+2$, this and the Littlewood--Paley inequality show that $T$ is bounded from $L^{d+2}$ to $L^{(d+2)/(d+1)}$.  
To prove \eqref{lp}, it is sufficient to show the  estimates 
\begin{align}
\label{l22}
\| P_jf\ast \sigma^{z} \|_2 &\le C (1+|\Im z|)^{3}  2^{-\frac{d}2 j}  \|f\|_2,  && \text{ if }  \re z=1/2, 
\\
\label{l10} \| P_jf\ast \sigma^{z} \|_\infty &\le C2^{j}\|f\|_1,   &&   \text{ if }  \re z= 0.
\end{align}
The first estimate \eqref{l22} is a straightforward consequence of  \eqref{osc-est0}  and Plancherel's theorem.  
The latter  \eqref{l10} is also clear, since $\|K_j\ast \sigma^{z}\|_\infty \lesssim 2^{j}$ if $ \re z= 0$.    Interpolation  along the analytic family  $f\mapsto   P_jf\ast \sigma^{z}$ 
yields  \eqref{lp}.

When $d=4$,   \eqref{osc-est0} and  Plancherel's theorem give  
\[ \| P_jf\ast \sigma^{z} \|_2 \le C (1+|\Im z|)^{3}  j 2^{-\frac{d}2 j}  \|f\|_2,  \text{ if }  \re z=1/2, \] 
instead of \eqref{l22} while \eqref{l10} remains valid. Following the same argument as above, we get the boundedness $(1/p, 1/q)\in \operatorname{int}\mathcal T$.

\subsection{Proof of Corollary \ref{cor:max}}   
We first show $L^p$ boundedness for $d>(d+1)/d$. 

Let $p>(d+1)/d$.   By a standard argument  it is sufficient to show 
\Be 
\label{max}
\| \sup_{t\in [1,2]}|P_jf\ast \sigma _t^{1/(d+1)}|\|_p\le C2^{-\epsilon j} \|f\|_p
\Ee
for some $\epsilon>0$  (e.g., see \cite{schlag}).  
The maximal operator is not  linear but sublinear, so  the interpolation along an analytic family is not directly applicable.  However, this can be  overcome 
by linearizing  the maximal operator. 

Let $x\mapsto \tau(x)$ be a measurable function from $\mathbb R^d$ to $[1,2]$.   
The estimate above is equivalent to 
\[ \| P_jf\ast \sigma _{\tau(\cdot)}^{1/(d+1)} \|_p\le C2^{-\epsilon j} \|f\|_p\]
with a constant $C$ independent of $\tau$, for  $p>(d+1)/d$. We consider  an analytic family of operators $f\mapsto   P_jf\ast \sigma _{\tau(\cdot)}^{z}$. 
By interpolation along the analytic family,  the desired estimate \eqref{max} follows if we prove  
 the  estimates
\begin{align*}
\| P_jf\ast \sigma _{\tau(\cdot)}^{z} \|_2 &\le C (1+|\Im z|)^{3}  j 2^{-\frac{d-1}2 j}  \|f\|_2,  &&  \re z=1/2, 
\\
\| P_jf\ast \sigma _{\tau(\cdot)}^{z}\|_1 &\le C2^{j} \|f\|_1,   &&     \re z= 0.
\end{align*}
Obviously,  these estimates also follow from the corresponding maximal estimates 
\begin{align*}
\| \sup_{t\in [1,2]}|P_jf\ast \sigma _{t}^{z} | \|_2 &\le C (1+|\Im z|)^{3}  j 2^{-\frac{d-1}2 j}  \|f\|_2,  &&   \re z=1/2, 
\\
\| \sup_{t\in [1,2]}|P_jf\ast \sigma _{t}^{z} \|_1 &\le C2^{j} \|f\|_1,   &&    \re z= 0.
\end{align*}
The first estimate is a consequence of the estimate \eqref{osc-est0}, Plancherel's theorem,  and the Sobolev imbedding. The second   is immediate  since $\sup_{t\in [1,2]} |K_j\ast \sigma _{t}^{z}(y)|\lesssim 2^j (1+|y|)^{-N}$.

\subsubsection*{Optimality of the damping order in Corollary \ref{cor:max}}\label{appendix22}    
Let $p\ge (d+1)/d$. 
We prove a slightly general statement that $\mathcal M^\rho$ generally  fails to be bounded on $L^q$ for all $q>p$,  if $\rho<1/(dp)$.  
Our example here is  a modification of  the one in \cite{Marletta}. 

 Let $\mathcal H$ be given by the graph of  $\phi(x)=1+|x|^{2m}$ for a positive integer $m$, which is clearly convex and analytic.  Let $q>p$, and  $\zeta_0$ be a nonnegative  bump function supported near the origin. We set 
\[ f(y',y_{d+1})=\zeta_0(y)|y_{d+1}|^{-1/q},  \quad y=(y', y_{d+1})\in \mathbb R^d\times \mathbb R,\] so 
$f\in L^p$.  
A computation shows that  $\H\phi(x)\sim |x|^{(2m-2)d}$. Hence, for a sufficiently small $\epsilon>0$, $|\mathcal M^\rho f(y)|\gtrsim\int_{B_{\epsilon}} f(y'-y_{d+1} x,  y_{d+1}|x|^{2m} )  |x|^{(2m-2)d\rho} dx$ provided that $|y|\le \epsilon$.  It follows that 
\[ 
    |\mathcal M^\rho f(y)|\gtrsim |y_{d+1}|^{-\frac 1q}\int_{B_{\epsilon}} |x|^{-\frac {2m}q+(2m-2)d\rho}dx
\]
for $|y|\le \epsilon$. 
If  $-\frac {2m}q+(2m-2)d\rho\le -d$,   $|\mathcal M^\rho f(y)|=\infty$ for $|y|\le \epsilon$.  
Thus,  in order that $\mathcal M^\rho$ is bounded  for all $q>p$, we must have  
\[ \rho\ge \frac{1}{(2m-2)d}\Big(\frac {2m}p-d\Big) \] for any integer $m\ge 1$. However, if $\rho< 1/(dp)$, by  taking $m$ large, the inequality above cannot be true. Consequently, $\mathcal M^\rho$ fails to be bounded on $L^q$  for all $q>p$. 

\section*{Acknowledgement} 
This work was supported by the National Research Foundation (Republic of Korea) grant RS-2024-00342160 (Lee) and the KIAS Individual Grant SP089101 (Oh).

\bibliographystyle{plain}

\end{document}